\documentclass[twocolumn]{autart}    

\usepackage{graphicx}
\usepackage[dvipsnames]{xcolor}
\graphicspath{{./figures/}}

\usepackage{amsmath,amssymb,amsfonts,enumitem,mathtools}

\usepackage{mathrsfs}
\usepackage{amsfonts,bbm}
\usepackage{dsfont}
\usepackage{epstopdf}
\usepackage{subfigure}
\usepackage[acronym]{glossaries}
\usepackage{balance}
\usepackage{url}

\newtheorem{theorem}{Theorem}
\newtheorem{definition}{Definition}
\newtheorem{proposition}{Proposition}
\newtheorem{lemma}{Lemma}
\newtheorem{corollary}{Corollary}
\newtheorem{example}{Example}
\newtheorem{remark}{Remark}

\newtheorem{standing}{Standing Assumption}

\newcommand{\R}{\mathbb{R}}

\newcommand{\D}{\mathbb{D}}

\newcommand{\N}{\mathbb{N}}

\newcommand{\mc}[1]{\mathcal{#1}}

\newcommand{\conv}{\mathrm{conv}}

\newcommand{\proj}{\mathrm{proj}}
\newcommand{\diag}{\mathrm{diag}}

\newcommand{\bs}{\boldsymbol}
\newcommand{\bsone}{\boldsymbol{1}}

\newcommand{\blue}{\textcolor{black}}

\newcommand\oprocendsymbol{\hbox{$\square$}}
\newcommand\oprocend{\relax\ifmmode\else\unskip\hfill\fi\oprocendsymbol}

\makeglossaries
\newacronym{wQLF}{wQLF}{weak Quadratic Lyapunov function}
\newacronym{wCLF}{wCLF}{weak Convex Lyapunov function}
\newacronym{wPLF}{wPLF}{weak Polyhedral Lyapunov function}
\newacronym{wLF}{wLF}{weak Lyapunov function}
\newacronym{EAS}{EAS}{Euler Auxiliary System}
\newacronym{LTI}{LTI}{linear time-invariant}
\newacronym{LTV}{LTV}{linear time-varying}
\newacronym{CT}{CT}{continuous-time}
\newacronym{DT}{DT}{discrete-time}
\newacronym[firstplural=Linear Matrix Inequalities (LMIs)]{LMI}{LMI}{Linear Matrix Inequality}
\newacronym{KSP}{KSP}{Kernel Sharing Property}

\begin{document}

\begin{frontmatter}

\title{Convergence in uncertain linear systems}

\author[delft]{Filippo Fabiani},
\author[tue]{Giuseppe Belgioioso},
\author[udine]{Franco Blanchini},
\author[milano]{Patrizio Colaneri},
\author[delft]{Sergio Grammatico}

\address[delft]{Delft Center for Systems and Control, TU Delft, The Netherlands}
\address[tue]{Control Systems group, TU Eindhoven, The Netherlands}
\address[udine]{Dipartimento di Scienze Matematiche, Informatiche e Fisiche, University of Udine, Italy}
\address[milano]{Dipartimento di Elettronica,  Informazione e Bioingegneria, Politecnico di Milano, and CNR-IEIIT, Italy}
\thanks{\emph{Email addresses:} \path{f.fabiani@tudelft.nl}~(F.~Fabiani), \path{g.belgioioso@tue.nl}~(G.~Belgioioso), \path{franco.blanchini@uniud.it}~(F.~Blanchini), \path{patrizio.colaneri@polimi.it}~(P.~Colaneri), \path{s.grammatico@tudelft.nl}~(S.~Grammatico).
	This work was partially supported by NWO under research projects OMEGA (613.001.702) and P2P-TALES (647.003.003), and by the ERC under research project COSMOS (802348).}

\begin{abstract}
State convergence is essential in several scientific areas, e.g. multi-agent consensus/disagreement, distributed optimization, monotone game theory, multi-agent learning over time-varying networks. \blue{This paper is the first on} state convergence in both continuous- and discrete-time linear systems affected by polytopic uncertainty. First, we characterize state convergence in linear time invariant systems via equivalent necessary and sufficient conditions.
In the presence of uncertainty, we complement the canonical definition of (weak) convergence with a stronger notion of convergence, which requires the existence of a common kernel among the generator matrices of the difference/differential inclusion (strong convergence). We investigate under which conditions the two definitions are equivalent. Then, we characterize weak and strong convergence by means of Lyapunov and LaSalle arguments, (linear) matrix inequalities and separability of the eigenvalues of the generator matrices. We also show that, unlike asymptotic stability, state convergence lacks of duality.
\end{abstract}

\end{frontmatter}


\section{Introduction}
Hand in hand with stability, state convergence of dynamical, possibly uncertain, systems represents a fundamental problem in system theory. However, while stability and asymptotic stability have been intensively studied in the system-and-control community, state convergence has received little attention, especially for uncertain systems. Unlike asymptotic convergence, where the state of a certain system is supposed to converge to a known (desired) state, with state convergence, we mean convergence of a system to some state, which in general is unknown a-priori. For this reason, we address the state convergence problem of continuous- and discrete-time linear systems, both time invariant and affected by polytopic uncertainty, by means of spectral and geometrical analysis, Lyapunov and LaSalle theories. \blue{To the best of our knowledge, we are the first to study convergence in uncertain linear systems.}
The problem of convergence to a constant, \emph{a-priori unknown} equilibrium state is ubiquitous and spans, for instance, from distributed consensus problems to multi-agent optimization and games over networks, positive system dynamics and tuning of plants with unknown input-output map. Two motivating problems are presented next. Further applications
are discussed in the example section (\S 8).

\textit{Motivating applications}:
We consider polar opinion dynamics in social networks \cite{Ghaderi2014,Amelkin2017,Tian2018}. Given a group of $N$ agents indexed by the set $\mc{I} \coloneqq \{1, \ldots, N\}$ and connected in a directed social network, we refer to the model proposed in \cite{Amelkin2017}, where the collective opinion vector $[x_1; \ldots; x_N] \eqqcolon x \in [-1,1]^N$, evolves according to
\begin{equation}\label{eq:opiniondynamics}
	\dot{x}(t) = - A(x(t)) L x(t).
\end{equation}
Here, $A(x(t)) \in \diag\left([0, 1]^N\right)$, for all $t \geq 0$, is the diagonal matrix that characterizes the susceptibility to persuasion of each agent, while $L \in \R^{N \times N}$ is the Laplacian matrix of the graph. The fundamental question in these models is whether or not the state $x(t)$ converges to some  \textit{a-priori unknown} state $\bar{x}$. In that case, we have consensus \blue{as} $A(\bar{x}) L \bar{x} = \bs{0}$, \blue{therefore} $\bar{x}$ is a zero of the mapping $x \mapsto A(x) L \blue{x}$.
Since $A(x(t)) \in \diag([0, 1]^N)$, $A(x) \in \conv\{A_1, \ldots, A_m\}$ for all $x \in [-1,1]^N$ and given matrices $A_i \in \R^{N \times N}$, $i \in \mc{M} \coloneqq \{1, \ldots, m\}$, then there always exist time-varying functions $\alpha_i(t, x(t)) \geq 0$, for all $i \in \mc{M}$ and $t \geq 0$, $\sum_{i \in \mc{M}} \alpha_i(t, x(t)) = 1$, such that $A(x(t)) = \sum_{i \in \mc{M}} \alpha_i(t, x(t)) A_i$. In \cite{Amelkin2017}, $A(x)$ is chosen equal to $\left(I - \diag(x)^2\right)$, $\tfrac{1}{2}\left(I - \diag(x)\right)$ or $\diag(x)^2$. Therefore, by denoting with $e_i$ the $i$-th vector of the canonical basis, in all the three cases, we have
$$A(x) \in \conv\left\{
\left[e_1,  \bs{0},  \ldots,  \bs{0}\right], \, \ldots,\, \left[\bs{0},  \ldots,   \bs{0},  e_N\right]\right\}.$$
Thus, the system in \eqref{eq:opiniondynamics} can be equivalently rewritten as
\begin{equation}\label{eq:op_dyn_poly}
\dot{x}(t) = - \left( \textstyle \sum_{i \in \mc{M}} \alpha_i(t, x(t)) A_i L \right) x(t).
\end{equation}
Clearly, the state solution to \eqref{eq:op_dyn_poly} belongs to the set of solutions to the linear differential inclusion $\dot{x} \in \conv\left(\{A_i L\}_{i \in \mc{M}}\right) x$, whose convergence implies convergence of the original nonlinear dynamics in \eqref{eq:opiniondynamics}. 

A second motivating application of state convergence is multi-agent learning over time-varying networks  
\cite{vamvoudakis2012multi,meng2015robust,meng2016learning,meng2017convergence}.
For example, in multi-agent games, the simplest dynamics to learn a Nash equilibrium are the projected pseudo-gradient dynamics, i.e., for each agent $i \in \mc{I}$,
\begin{equation} \label{eq:FB}
x_i(k+1) = \proj_{\Omega_i}\big( x_i(k) - \epsilon \nabla_{x_i} J_i(x_i(k), \bs{x}_{-i}(k)) \big),
\end{equation}
where $\Omega_i$ is a local constraint set, $\epsilon > 0$ is a step size and $J_i$ is the local cost function, which depends on the local decision variable $x_i$ (first argument) and on the decision variables of the other agents $\bs{x}_{-i}$ (second argument).
For simplicity, let $x_i$ be of dimension $1$. Let $J_i$ be quadratic in $x_i$ and affine in $\bs{x}_{-i}$, hence $\nabla_{x_i} J_i(x_i, \bs{x}_{-i}) = \Phi_i^\top \bs{x} + \phi_i$, for some $\Phi_i \in \R^{N}$ and $\phi_i \in \R$, and $\Omega_i = [a_i, b_i] \supset \{0\}$, hence the projection is a switched affine operator, i.e., $\proj_{\Omega_i}(v) \in \left\{ a_i, v, b_i  \right\} $.
Therefore, the dynamics in \eqref{eq:FB} are contained in switched affine dynamics $x_i^+ = \proj_{\Omega_i}\big( x_i - \epsilon \nabla_{x_i} J_i(x_i, \bs{x}_{-i}) \big) \in \left\{  a_i,  x_i - \epsilon \, \Phi_i^\top \bs{x} - \epsilon \, \phi_i, b_i \right\}$, and in turn in the uncertain affine dynamics
\begin{align*}
y_i(k)   & = \, x_i(k) - \epsilon \, \Phi_i^\top \bs{x}(k) - \epsilon \, \phi_i  \\
x_i(k+1) & =
\left\{ 
\begin{array}{ll}
y_i(k) & \textup{ if } a_i < y_i(k) < b_i \\
\displaystyle (b/{y_i(k)}) \, y_i(k) & \textup{ if } y_i(k) \geq b_i \\
\displaystyle (a/{y_i(k)}) \, y_i(k) & \textup{ if } y_i(k) \leq a_i  \\
\end{array}
\right. \\ 
 & \in [0,1] \left( x_i(k) - \epsilon \, \Phi_i^\top \bs{x}(k) - \epsilon \, \phi_i \right)\,,
\end{align*}

\vspace{-0.25cm}

since $b/{y_i(k)}, \, a/{y_i(k)} \in [0,1]$ in their respective intervals.
Thus, for the collective dynamics, we have the following affine difference inclusion:
\begin{align} \label{eq:projPG}
\bs{x}(k+1) & = \proj_{\Omega}  \nonumber
( \bs{x}(k)- \epsilon \, F (\bs{x}(k)) \\
       & \in [0,1]^N \left( (I - \epsilon \, \Phi^\top) \, \bs{x}(k) - \epsilon \, \phi \right)\,,
\end{align}
where $\Omega \coloneqq \times_{i=1}^N \Omega_i$ and $F(\bs{x}) \coloneqq \times_{i=1}^N \nabla_{x_i} J_i(x_i, \bs{x}_{-i})$ is the so-called pseudo-gradient mapping.
Clearly, convergence of the affine difference inclusion in \eqref{eq:projPG} implies convergence of the original nonlinear dynamics in \eqref{eq:FB}. In this paper, we study uncertain \textit{linear} systems, since the \textit{affine} case can be cast into the linear one via additional auxiliary states, see \S $8.1$.

Essentially, the Nash equilibrium problem is the problem to ensure that the state $\bs{x}(t)$ of the system in \eqref{eq:projPG} converges to some \textit{a-priory unknown} state $\bs{x}^*$, which happens to be a Nash equilibrium due to the structure of the dynamics.
Similarly, primal-dual projected pseudo-gradient dynamics can be designed for computing a generalized Nash equilibrium in games with coupling constraints \cite{BelGra17,yi:pavel:19,BelGra18}.
While multi-agent optimization and game equilibrium dynamics are typically nonlinear, the analysis of the linear case provides necessary conditions for convergence \cite{BFBG18}, \blue{e.g.} potential certificates of non-convergence, see \cite{grammatico18} for an example of non-convergent linear time-varying primal-dual dynamics.

\textit{Contribution}:

\vspace{-3mm}

\begin{itemize}
	\item \blue{Preliminarily, we start with} necessary and sufficient conditions for state convergence of linear time invariant systems. Specifically, we link state convergence with the existence of a (weak) Lyapunov function, the separability of the eigenvalues, the stability of an auxiliary system and (linear) matrix inequalities (\S3 - for ease of readability, the proofs of this section are reported in Appendix~\ref{sec:app_linear});
	\item We introduce the notions of weak and strong convergence for uncertain \blue{linear} systems. We show that a sufficient condition, i.e., the kernel sharing property among the generator matrices of the difference/differential inclusion, implies that the two definitions are equivalent (\S4);
	\item \blue{We show that the existence of a quadratic or polyhedral Lyapunov function is not sufficient for strong convergence: the existence of a suitable common decomposition of the generating matrices is required (\S5)};
	\item We associate sufficient (linear) matrix inequalities to weak and strong convergence (\S5-6);
	\item \blue{We investigate weak convergence via LaSalle arguments and define the concept of weak kernel (\S6);}
	\item We show the lack of duality of state convergence in uncertain \blue{linear} systems. \blue{Nevertheless}, we show that the existence of a quadratic Lyapunov function guarantees duality of strong convergence (\S7).
\end{itemize}

\textit{Notation}: 
\blue{$\N,$
$\R$ and $\R_{\geq 0}$ denote the sets of natural, real and non-negative real numbers, respectively;
$\mathrm{bdry}(S)$ denotes the boundary of a set $S$. $\textrm{ker}(A)$
denotes the kernel of matrix $A$;
$\mu(A)$ denotes its Lozinski measure, i.e., $\mu(A)=\lim_{h\to 0} \frac{\|I+hA\|-1}{h}$;
$P \succ (\succcurlyeq) \, 0$ denotes that $P$ is a positive (semi-)definite symmetric matrix.
$\bs{0}$ ($\bs{1}$) denotes vectors with elements all equal to $0$ ($1$).
}

\section{Mathematical background}

\blue{Let us} consider a \gls{DT} \gls{LTI} system
\begin{equation}\label{eq:LTI_DT}
	x(k+1) = A \, x(k),
\end{equation}

where $x \in \R^n$ and $A \in \R^{n \times n}$, $k \in \N$. For \gls{LTI} systems, we are interested in global convergence of the state $x$ to some vector $\bar{x}$, which may depend on the initial conditions $x(0)$, accordingly with the following definition.

\begin{definition}\label{def:convergence}
	\textup{(\textbf{Convergence})}\\
	The system in \eqref{eq:LTI_DT} is convergent if, for all $x(0) \in \R^n$, there exists $\bar{x} \in \R^n$ such that the solution $x(k)$ to \eqref{eq:LTI_DT} satisfies $\lim_{k \rightarrow \infty} \left\| x(k) - \bar{x} \right\| = 0$.
	\hfill $\square$
\end{definition}

The convergence of \gls{LTI} systems is closely related with the location of the eigenvalues of the matrix $A$ and their algebraic/geometric multiplicity. Thus, \blue{let us} recall the notion of semi-simple eigenvalue \cite{BCM09} and of (weak) Lyapunov function.

\begin{definition}
	\textup{(\textbf{(Semi-) Simple eigenvalue})}\\
	An eigenvalue is \textit{semi-simple} if it has equal algebraic and geometric multiplicity. An eigenvalue is \textit{simple} if it has algebraic and geometric multiplicities both equal to $1$.
	\hfill $\square$
\end{definition}


\begin{definition}\label{def:wlf}
	\textup{(\textbf{Weak Lyapunov function})}\\
	A positive definite function $V: \R^n \rightarrow \R_{\geq 0}$ is a \gls{wLF} for the system in \eqref{eq:LTI_DT} if $V(A x) \leq V(x)$, for all $x \in \R^n$.
	\hfill $\square$
\end{definition}

In particular, to link convergence results and stability of dynamical systems, we will consider several classes of \glspl{wLF}, as \glspl{wQLF}, \glspl{wPLF} and, more generally, \glspl{wCLF}.

The same definitions of convergence and stability apply to the \gls{CT} \gls{LTI} system
 
\vspace{-0.5cm}  
\begin{equation}\label{eq:LTI_CT}
	\dot{x}(t) = A \, x(t),
\end{equation}
\vspace{-0.5cm} 

by substituting $k \in \N$ with $t \geq 0$. The definition of \gls{wLF} in Definition~\ref{def:wlf}, instead, reads as $\dot{V}(x) \leq 0$ for all $x \in \R^{n}$. 

\section{Convergence in LTI systems}\label{sec:LTI}

\subsection{Discrete-time systems}
\blue{We start with} some equivalence results that link the convergence of \eqref{eq:LTI_DT} with the \blue{(at least marginal)} stability of an auxiliary system, \blue{(linear) matrix inequality (LMI)} conditions and existence of a \gls{wLF}. The proofs are reported in the appendix.

\begin{lemma}
	\label{lemma:convergent-stable}
	The system in \eqref{eq:LTI_DT} is convergent if and only if there exists $\eta \in (0,1)$ such that the system 
	\begin{equation}\label{eq:LTI_DT_aux}
		x(k+1) = \left(  \frac{1}{\eta} A - \frac{1-\eta}{\eta}I \right) x(k) \eqqcolon A^{\textup{dt}}_\eta x(k),
	\end{equation}
	 is \blue{(marginally)} stable.
{\hfill $\square$}
\end{lemma}

\begin{proposition}\label{th:LTI_DT}
	\textup{(\textbf{Convergence in \gls{DT} LTI systems})} \smallskip \\
	The following statements are equivalent:
	\begin{enumerate}[label=(\alph*)]
		\item The system in \eqref{eq:LTI_DT} is convergent;
		\smallskip
		\item $\exists \, \eta \in (0,1)$ such that  the system in \eqref{eq:LTI_DT_aux} is \blue{(marginally)} stable;
		\smallskip
		\item there exists a \gls{wPLF} for the system in \eqref{eq:LTI_DT_aux};
		\smallskip
		\item there exists an invertible matrix $T \in \R^{n \times n}$ such that 
		$$T^{-1} A T = 
		\left[ \begin{array}{ccc}
		A^{\textup{as}} & & 0 \\ 
		A^{\textup{r}} & & I 
		\end{array} \right],
		$$
		for some Schur matrix $A^{\textup{as}}$;

\blue{ 
\item $\exists \, \eta \in (0,1)$ and $P \succ 0$ such that:
		$$
		\eta (A^\top PA-P)+(1-\eta)(A^\top -I)P (A-I) \preccurlyeq 0 \,;
		$$
	\item 
	$\exists Q \succ 0$ and $P \succcurlyeq 0$, with ${\rm rank}(P) = n - {\rm dim}({\rm ker}(A-I))$,  such that:
		$$
		A^\top PA-P+(A^\top -I)Q (A-I) \preccurlyeq 0 \,;
		$$
}
{\hfill $\square$}
	\end{enumerate}
\end{proposition}

\blue{We note that if the matrix inequality $(e)$ in Proposition \ref{th:LTI_DT} holds for some $(\eta, P)$, then it also holds for $(\eta_1, P)$, where $\eta < \eta_1 < 1$. This is because the function $\eta \mapsto -(1-\eta)/\eta$ is increasing for $\eta>0$. Thus, we can fix $\eta$ arbitrarily close to $1$ and solve the corresponding LMI for $P \succ 0$.}

\begin{remark}
The system in \eqref{eq:LTI_DT} is convergent if and only if the mapping $x \mapsto A x$ is ``averaged'' \cite{BFBG18}. 
{\hfill$\square$}
\end{remark}


\blue{
We conclude the subsection by noticing that the existence of a wPLF does not imply convergence. For example, consider the system $x(k+1) = \left[ \begin{smallmatrix}
0 & -1  \\
1 &  \phantom{-}0
\end{smallmatrix} \right] x(k)$, where the state does not converge, yet $\| x(k) \|_\infty$ is constant.
}


\subsection{Continuous-time systems}
To characterize the convergence of the \gls{CT} \gls{LTI} in \eqref{eq:LTI_CT}, let us \blue{also} introduce the following \gls{DT} auxiliary system.

\begin{definition}\label{def:}
	\textup{(\textbf{Euler Auxiliary System})}\\
	Given $\tau > 0$, the \gls{DT} system
	\begin{equation}\label{eq:EAS}
	x(k+1) = \left( I + \tau A\right) x(k) \eqqcolon A_\tau x(k)
	\end{equation}
	is the \gls{EAS} of the \gls{CT} system in \eqref{eq:LTI_CT}.
	{\hfill $\square$}
\end{definition}

\blue{We then have} the following results.

\begin{lemma} \label{lemma:convergent-stable_tc}
	The system in \eqref{eq:LTI_CT} is convergent if and only if there exists $\epsilon > 0$ such that the system 
	\begin{equation}\label{eq:LTI_CT_aux}
		\dot{x}(t) = A \, (I + \epsilon A)^{-1}  \, x(t) \eqqcolon A^{\textup{ct}}_\epsilon \, x(t),
	\end{equation}
	is \blue{(marginally)} stable.
{\hfill $\square$}
\end{lemma}

\begin{proposition}\label{th:LTI_CT}
	\textup{(\textbf{Convergence in \gls{CT} LTI systems})} \smallskip \\
	The following statements are equivalent:
	\begin{enumerate}[label=(\alph*)]
		\item The system in \eqref{eq:LTI_CT} is convergent;
		\smallskip
		\item  $\exists \, \epsilon > 0$ such that the system in \eqref{eq:LTI_CT_aux} is \blue{(marginally)} stable;
		\smallskip
		\item there exists a \gls{wPLF} for the system in \eqref{eq:LTI_CT};
		\smallskip
		\item there exists an invertible matrix $T \in \R^{n \times n}$ such that 
		$$
		T^{-1} A T = \left [
		\begin{array}{ccc} 
		A^{\textup{as}} & & 0 \\ 
		A^{\textup{r}}  & & 0 
		\end{array} \right ],
		$$
		for some Hurwitz matrix $A^{\textup{as}}$;
		\smallskip
		\item $\exists \, \tau>0$ such that the \gls{EAS} in \eqref{eq:EAS} converges;
		\smallskip
		\item $\exists \, \epsilon > 0$ and $P \succ 0$ such that:	
		$$
		A^\top P + PA + \epsilon \, A^\top P A \preccurlyeq 0 \,;
		$$
		\smallskip
 		\blue{
\item  $\exists Q \succ 0$ and $P\succcurlyeq 0$, with ${\rm rank}(P)=n-{\rm dim}({\rm ker}(A))$,  such that:
		$$
		A^\top P + PA + A^\top Q A \preccurlyeq 0 \,.
			\vspace{-0.75cm}
		$$
}
		\hfill $\square$
	\end{enumerate}
\end{proposition}

\blue{
Similarly to Proposition \ref{th:LTI_DT} $(e)$, we can fix an arbitrarily small $\epsilon > 0$ and solve the corresponding LMI $(f)$ in Proposition \ref{th:LTI_CT} for $P \succ 0$.
}

\section{Uncertain linear systems}
\label{sec:LU}
For difference or differential linear inclusions \cite{Aub} the definition of convergence requires some care. Specifically, we shall distinguish between weak and strong convergence. 
\subsection{Weak and strong convergence}
We consider uncertain \gls{DT} linear systems of the form
\begin{equation}\label{eq:polytopicdiffincl}
x(k+1) = A(w(k)) x(k),
\end{equation}
where $A(w(k))$ satisfies the following assumption.
\begin{standing}
	\textup{(\textbf{Polytopic uncertainty})}\\
	$$A(w) \coloneqq \textstyle \sum_{i \in \mc{M}} A_i w_i,$$
	with $\mc{M} \coloneqq \{1, 2, \ldots, M \}$ and $w \in \mc{W}$, defined as
	$$
	\mc{W} \coloneqq \left\{ w \in \R^M \, \Big| \, \sum_{i \in \mc{M}} w_i = 1, \; w_i \geq 0 \; \forall i \in \mc{M} \right\}.
	\vspace{-0.75cm}
	$$ 
	{\hfill$\square$}
\end{standing}
In \gls{CT}, we consider the differential inclusion of the form
\begin{equation}\label{eq:polytopicdiffinclcont}
\dot{x}(t) = A(w(t)) x(t),
\end{equation}
with the same polytopic uncertainty structure.
As for \gls{LTI} systems, we investigate whether $x(k)$ converges to some $\bar{x}$, which in general depends on $w(k)$. For instance, if $ A(w(k)) = \left[ \begin{smallmatrix} a_{1,1}(k) & \; & 0 \\ 1 & \; & 1 \end{smallmatrix}\right]$, 
with $a_{1,1}(k) \in \{1/2, \, 3/4\}$ for all $k \in \N$, then $x_1(k) \rightarrow 0$ and $x_2(k) \rightarrow \overline{x}_2 = x_2(0) + \sum_{k\geq 0} x_1(k)$, which is finite but depends on $\left(a_{1,1}(k)\right)_{k \in \N}$.

In view of this example, we give two different definitions of convergence.
\begin{definition}
	\textup{(\textbf{Weak convergence})}\\
	The difference inclusion in \eqref{eq:polytopicdiffincl} (respectively, differential inclusion in \eqref{eq:polytopicdiffinclcont}) is weakly convergent if, for all sequences  $w(k) \in \mc{W}$ (resp., for all $w(t) \in \mc{W}$) and initial conditions $x(0) \in \R^{n}$, there exists a vector $\bar x \in \R^n$ such that $\underset{ k \rightarrow \infty }{\lim} \left\| x(k) -   \bar x \right\| = 0$ $\left(\underset{t \rightarrow \infty}{\lim} \left\| x(t) -   \bar x \right\| = 0\right)$.
	{\hfill $\square$}
\end{definition}

Next, we introduce a stronger notion of convergence, i.e., convergence to the common kernel of the matrices $\{A_i-I\}_{i \in \mc{M}}$ in DT \eqref{eq:polytopicdiffincl}, or $\{A_i\}_{i \in \mc{M}}$ in CT \eqref{eq:polytopicdiffinclcont}.

\begin{definition}
	\textup{(\textbf{Strong convergence})}\\
	The difference inclusion in \eqref{eq:polytopicdiffincl} (resp., differential inclusion in \eqref{eq:polytopicdiffinclcont}) is strongly convergent if, for all sequences $w(k) \in \mc{W}$ (resp., for all $w(t) \in \mc{W}$) and initial conditions $x(0) \in \R^{n}$, there exists a vector $\bar x \in \bar{\mathcal{X}} \subseteq \R^n$, where 
	$$
	\bar{\mathcal{X}} \coloneqq \bigcap_{i \in \mc{M}}~\ker\left(A_i - I\right), \quad \left ( \text{resp., } \bar{\mathcal{X}} \coloneqq \bigcap_{i \in \mc{M}}~\ker(A_i) \right)
	$$ 
	such that $\displaystyle \lim_{k \rightarrow \infty} \left\| x(k) - \bar x \right\| = 0$ $\left( \displaystyle\lim_{ t \rightarrow \infty } \left\| x(t) -   \bar x \right\| = 0\right)$.
	{\hfill $\square$}
\end{definition}

Therefore, to have strong convergence, the limit vector must be a steady state.
In other words, if we initialize the state $x(0) \in \bar{\mathcal{X}}$ then $x(k) = x(0)$ for all $k > 0$ and for all possible sequences $w(k) \in \mc{W}$. 
On the other hand, this is not ensured with a limit vector in the case of weak convergence. Moreover, as stressed in the following example, while strong convergence implies weak convergence, the converse does not necessarily hold. 

\begin{example}
	The scalar difference inclusion $x(k+1) \in \{1/2, \; 1\} \, x(k)$ is weakly convergent (because the sequence $\{ x(k) \}_{k \in \N}$ is non-increasing over $k \geq 0$) but not strongly convergent.
	{\hfill$\square$}
\end{example}

\subsection{Kernel sharing property: when weak convergence implies strong convergence}
We now investigate under which conditions the weak convergence of an uncertain system of the form \eqref{eq:polytopicdiffincl} (or \eqref{eq:polytopicdiffinclcont}) implies strong convergence. 

\begin{proposition}\label{prop:strongimpliessharing}
	If the difference inclusion in \eqref{eq:polytopicdiffincl} (resp., differential inclusion in \eqref{eq:polytopicdiffinclcont})
	is strongly  convergent, then the matrices $\left\{A_i - I\right\}_{i \in \mc{M}}$ (resp.,  $\left\{A_i\right\}_{i \in \mc{M}}$), have exactly the same kernel:
	\begin{equation}\label{commonker}
	\begin{aligned}
	&\,\ker\left(A_i - I\right) = \bigcap_{i \in \mc{M}} \ker\left(A_i - I\right), \quad \forall i \in \mc{M}\\
	&\left(\text{resp., }  \ker(A_i)  = \bigcap_{i \in \mc{M}} \ker(A_i), \quad \forall i \in \mc{M}\right).
	\end{aligned}
\vspace{-0.75cm}
	\end{equation}
{\hfill $\square$}
\end{proposition}
\begin{pf}
	We prove the \gls{DT} case by contradiction (the proof in \gls{CT} is analogous). Let us assume that there exists some $j \in \mc{M}$ such that we have strict inclusion, i.e., $\bar{\mc{X}} \subset  \ker\left(A_j - I\right)$. Then, for some $x(0) \in \ker\left(A_j - I\right) \setminus\bar{\mc{X}}$, and for $A(w(k)) = A_j$ for all $k \in \N$, we have the constant solution $x(k) = x(0) \notin \bar{\mc{X}}$.	
	\hfill$\blacksquare$
\end{pf}

Along this idea, we link weak and strong convergence via the \gls{KSP} introduced next.

\begin{definition}\textup{(\textbf{Kernel Sharing Property})}\\
	The family of matrices $\{A_i\}_{i \in \mc{M}}$ has the Kernel Sharing Property if \eqref{commonker} holds true.
	\hfill$\square$
\end{definition}

Note that the \gls{KSP} is equivalent to claiming that $ \ker \left(A(w)\right ) = \mc{K} $ (resp., $\ker \left(A(w)-I \right ) = \mc{K}$ ) is invariant with respect to $w \in \mc{W}$.
We now show that, under the \gls{KSP} assumption, if the system is weakly convergent, then it is also strongly convergent. Some preliminary results are reported first.  

\begin{lemma}\label{gotozero_dis}
	Let $x(k)$ be a solution to \eqref{eq:polytopicdiffincl}. If $x(k) \rightarrow \bar x$, then $\left(A(w(k)) - I\right)\bar x \rightarrow 0$.
\hfill$\square$
\end{lemma} 
\begin{pf}
	Let $z(k) \coloneqq x(k) -\bar x \, \blue{\rightarrow 0}$. \blue{By} \eqref{eq:polytopicdiffincl}, we have 
	$$
	z(k+1) = A(w(k)) z(k) + \left(A(w(k)) - I\right) \bar x.
	$$
	\blue{Then}, $z(k) \rightarrow 0$ \blue{implies that} $\left(A(w(k)) - I\right) \bar x\rightarrow 0$.
	\hfill$\blacksquare$
\end{pf}

The \gls{CT} counterpart is not immediate, due to some technical problems, as stressed by the following example.
\begin{example}
	The scalar differential inclusion $\dot x(t) \in \{0, \; -1\} \, x(t)$ is weakly convergent (because $x(t)$ is non-increasing over $t \geq 0$). Now, let us consider the case that $A$ switches to $-1$ with spikes of decreasing length $\Delta_h < 1$, i.e., let $A(w(t))=-1$ when $t \in [h, h+\Delta_h]$, and $A(w(t))=0$ otherwise.
	By considering an interval size such that $\sum_{h \geq 0}\Delta_h =\ln(2)$, we observe that $x(t)$ tends asymptotically to $x(0)/2$, which is not in the common kernel of the $A_i$ matrices.
	{\hfill$\square$}
\end{example}

Thus, the natural extension of Lemma~\ref{gotozero_dis} to the \gls{CT} case does not hold, i.e., $x(t) \rightarrow \bar x$ does not imply that $A(w(t)) \bar x \rightarrow 0$. However, convergence does happen on average, as formalized next.

\begin{lemma}\label{gotozero_con}
	Let $x(t)$ be a solution to \eqref{eq:polytopicdiffinclcont}. If $x(t) \rightarrow \bar x$, then 
	$$
	\lim_{T \to \infty} \, \frac{1}{T} \int_0^T A(w(t)) \bar{x} \; dt = 0 \,.
	\vspace{-0.5cm}
	$$
\hfill$\square$
\end{lemma} 
	\vspace{-0.5cm}
\begin{pf} 
	Let $z(t) \coloneqq x(t) -\bar x \to 0$. Then, by \eqref{eq:polytopicdiffinclcont},
	$$
	\dot z(t) = A(w(t)) z(t) + A(w(t)) \bar x.
	$$
	By integrating and diving by $T > 0$, we obtain
	\begin{equation*}
	\begin{aligned}
	&\frac{1}{T}\int_0^T A(w(t)) \bar{x} \; dt = \frac{1}{T} \int_{0}^{T} \left (\dot{z}(t)-A(w(t))z(t) \right)  \; dt\\
	&= \frac{1}{T}[x(T) -x(0)] - \frac{1}{T}\int_0^T A(w(t))z(t) \;  dt \to 0.
	\end{aligned}
	\end{equation*}
	For $T \to \infty$, the first term converges to $0$ because $x(t)$ is bounded and, since $A(w(t)) z(t) \to 0$,
        the average of a function converging to $0$ converges to 0 as well. Thus, $\lim\limits_{T \to \infty} \frac{1}{T} \int_{0}^{T} A(w(t)) \bar{x} \; dt = 0$.
	\hfill$\blacksquare$
\end{pf}
\begin{lemma}\label{must_be_zero}
	Let $x(k)$ be a solution to \eqref{eq:polytopicdiffincl} (resp., $x(t)$ solution to \eqref{eq:polytopicdiffinclcont}). If $x(k) \rightarrow \bar x$ (resp., $x(t) \rightarrow \bar x$), then there exists some $\bar w \in \mc{W}$ such that $\left(A(\bar w) - I\right) \bar x=0$ (resp., $A(\bar w) \bar x = 0$).
	\hfill$\square$
\end{lemma}
	\vspace{-0.5cm}
\begin{pf}
	In the \gls{DT} case, by Lemma~\ref{gotozero_dis} we have:
	$$ 
	u(k) \coloneqq \sum_{i \in \mc{M}}^{} \left(A_i w_i(k) - I\right) \bar{x} = \sum_{i \in \mc{M}}   w_i(k)  \left(A_i - I \right)\bar x  \rightarrow 0.
	$$
	For all $k \in \N$, since $w(k) \in \mc{W}$, the vector $u(k)$ belongs to the convex hull of the vectors $\left\{(A_i  - I)\bar x\right\}_{i \in \mc{M}}$, which is closed and convex. Thus, in the limit for $k \to \infty$, $u(k)$ shall belong to the convex hull as well, i.e.,
	$\textstyle
	0 = \sum_{i \in \mc{M}}   \bar w_i   \left (A_i  - I\right) \bar x,
	$
	namely, $\left( A(\bar{w}) - I \right) \bar x =0$.
	
	In the \gls{CT} case, by Lemma~\ref{gotozero_con}, we have:
	\begin{equation*}
	\begin{aligned}
	&\frac{1}{T} \int_0^T \sum_{i \in \mc{M}}   w_i(t)  A_i  \bar x \; dt = 
	\left(\sum_{i \in \mc{M}} \frac{1}{T} \int_0^T  w_i(t) \; dt \right)  A_i  \bar x\\
	&= \sum_{i \in \mc{M}} \omega_i(T)  A_i  \bar x  \rightarrow 0,
	\end{aligned}
	\end{equation*}
	where $\omega_i(T) \coloneqq \frac{1}{T}\int_0^T  w_i(t) \; dt$ denotes the nonnegative average values of $\{w_i(t)\}_{i \in \mc{M}}$ for every $T > 0$.
	Then, to derive $A(\bar w) \bar x=0$ for some $\bar w$, we use the same argument on the limit used in the \gls{DT} case.
	\hfill$\blacksquare$
\end{pf}

	\vspace{-0.5cm}
	
We are now ready to exploit the previous lemmas to link weak and strong convergence under the \gls{KSP}.
\begin{proposition} \label{prop:sharingimpliesstrong}
	If the family of matrices $\{A_i\}_{i \in \mc{M}}$ has the Kernel Sharing Property, then weak convergence implies strong convergence.
\hfill$\square$
\end{proposition}
	\vspace{-0.5cm}
\begin{pf}
	In view of Lemma~\ref{must_be_zero}, the convergence (if it holds)
	shall be in the kernel of some $(A(\bar w)-I)$ (resp., $A(\bar w)$), since the kernel is common.
	\hfill$\blacksquare$
\end{pf}

\begin{remark}
	In view of the properties of \blue{the Laplacian matrix of a strongly connected graph,}
	the family of matrices $\{A_i L\}_{i \in \mc{M}}$ in \eqref{eq:op_dyn_poly}, \blue{the first motivating application}, has the \gls{KSP}, since $\bsone \in \ker(A_i L)$, for all $i \in \mc{M}$.
	{\hfill$\square$}
\end{remark}

The \gls{KSP}, which can be efficiently verified by means of linear algebra arguments, is crucial to investigate strong convergence. In fact, it follows from Propositions~\ref{prop:strongimpliessharing} and \ref{prop:sharingimpliesstrong} that if the matrices do not have a common kernel, there can not be strong convergence.
\section{Strong convergence}

\subsection{Separability of the eigenvalues}

We first characterize the strong convergence property of difference (differential) inclusions via separability of the eigenvalues of the generator matrices.
Essentially, if the \gls{KSP} holds, then the convergence analysis reduces to investigate the asymptotic stability of a subsystem.
\begin{theorem}\label{th:differentialinclusions_strongconv}
	The following statements are equivalent:
	\begin{enumerate}[label=(\alph*)]
		\item The difference inclusion in \eqref{eq:polytopicdiffincl} (resp., differential inclusion in \eqref{eq:polytopicdiffinclcont}) is strongly convergent;
		\smallskip
		\item there exists an invertible matrix $T \in \R^{n \times n}$ such that, for all $i \in \mc{M}$,
		\begin{equation}\label{eq:partition}
		 T^{-1} A_i T = 
		\left[ \begin{array}{ccc}
		A_i^{\textup{as}} & & 0 \\
		A_i^{\textup{r}} & & I_m
		\end{array} \right]\;
		 \left(\text{resp., }=\,
		\left[ \begin{array}{ccc} 
		A_i^{\textup{as}} & & 0 \\ 
		A_i^{\textup{r}} & & 0_m 
		\end{array} \right]
		\right),
		\end{equation}
		where $m \coloneqq \dim(\bar{\mc{X}})$ and the matrices $\{A_i^{\textup{as}}\}_{i \in \mc{M}}$ generate an asymptotically stable difference inclusion (resp., differential inclusion).
{\hfill $\square$}
	\end{enumerate}
\end{theorem}
\begin{pf}
	We give the proof for the \gls{DT} case, since the one in \gls{CT} is analogous. Let $T_2 \in \R^{n \times m}$ be the matrix generated by a basis of $\bar{\mc{X}}$, i.e., $\left(A_i - I\right) T_2 = 0$ for all $i \in \mc{M}$, and let $T_1 \in \R^{n \times (n-m)}$ be a complement of $T_2$ such that $T \coloneqq [T_1 \; T_2] \in \R^{n \times n}$ is invertible. Then, $T$ determines a similarity transformation that separate the eigenvalues of the matrices $\{A_i\}_{i \in \mc{M}}$ as follows: 
		$$
	T^{-1} A_i T=\left[\begin{array}{ccc}
	A_{i1} & & 0\\
	A_{i2} & & I_m
	\end{array}\right], \quad \text{ for all } i \in \mc{M}.
	$$ 
	Since the difference inclusion \eqref{eq:polytopicdiffincl} converges to $\bar{\mc{X}}$, this implies that the matrices $\{A_{i1}\}_{i \in \mc{M}}$ shall generate an asymptotically stable difference inclusion.
	{\hfill$\blacksquare$}
\end{pf}

\subsection{Lyapunov-like results}

Next, we characterize the strong convergence property for both difference and differential inclusions by means of \blue{a} Lyapunov analysis. First, we show that strong convergence implies the existence of a polyhedral Lyapunov function.
\begin{theorem}\label{th:strongconv-wplf_differentialinclusion}
	If the difference inclusion in \eqref{eq:polytopicdiffincl} (resp., differential inclusion in \eqref{eq:polytopicdiffinclcont}) is strongly convergent, then it admits a \gls{wPLF}.
	{\hfill$\square$}
\end{theorem}
	\vspace{-0.5cm}
\begin{pf}
	We give the proof for the \gls{CT} case (the one in \gls{DT} is analogous).
	The strong convergence of \eqref{eq:polytopicdiffinclcont} allows to partition each matrix $A_i$ as in \eqref{eq:partition}
	with $\{A_i^{\textup{as}}\}_{i \in \mc{M}}$ that generate an asymptotically stable differential inclusion. 
	By \cite[Prop.~7.39]{BM15}, it admits a PLF whose vertices are in $X \in \R^{(n-m) \times r}$ such that, for all $i \in \mc{M}$,
	$$
	A_i^{\textup{as}} X = X P_i,
	$$
	with $P_i \in \R^{r \times r}$ strictly column diagonally dominant, i.e., $\mu(P_i)< 0$, for all $i \in \mc{M}$. Then, let us consider the augmented system with $\beta > 0$:
	$$
	\left[\begin{array}{ccc}
	A_i^{\textup{as}} & & 0\\
	A_i^{\textup{r}} & & 0_m
	\end{array}\right] \,
	\overbrace{\left[\begin{array}{ccc}
		X & & 0\\
		0 & & \beta I_m
		\end{array}\right]}^{\eqqcolon X^{\textup{aug}}} = 
	\left[\begin{array}{ccc}
	X & & 0\\
	0 & & \beta I_m
	\end{array}\right] \,
	\overbrace{\left[\begin{array}{ccc}
		P_i & & 0\\
		\frac{1}{\beta} A_i^{\textup{r}} X & & 0_m
		\end{array}\right]}^{\eqqcolon P_i^{\textup{aug}}}.
	$$
	For $\beta$ sufficiently large, each matrix $P_i^{\textup{aug}}$ is such that $\mu(P_i^{\textup{aug}})\le 0$, so that  $X^{\textup{aug}}$ represents the matrix of vertices of a \gls{wPLF}.	
	{\hfill $\blacksquare$}
\end{pf}
	\vspace{-0.5cm}

While strong convergence implies the existence of a \gls{wPLF}, it does not imply the existence of a \gls{wQLF}. This follows by seeing asymptotic stability as a special case of strong convergence to the kernel $\{0\}$ (i.e., robust asymptotic stability). \blue{In fact}, there are \blue{uncertain linear}  systems that are asymptotically stable but do not admit a quadratic Lyapunov function \cite{BM15}. \blue{Conversely}, the existence of either a \gls{wPLF} or a \gls{wQLF} in general does not imply convergence, \blue{neither strong nor weak, as shown by the following examples}.

\blue{
\begin{example}
For the \gls{DT} system $x(k+1) \in \{-1,\; 1\} x(k)$, with $x(0) = 1$, $V(x) = x^2$ is a \gls{wQLF} and $\sqrt{V(x)} = |x|$ is a \gls{wPLF}. Whenever $A(w(k)) = 1$ for $k$ even and $A(w(k)) = -1$ for $k$ odd, the system does not converge.
\hfill $\square$
\end{example}
}

\blue{
\begin{example} \label{ex:QLFPLF_not_conv}
Consider a \gls{CT} system with $A(w(t)) \in \left\{ 
	\left[
	\begin{smallmatrix}
	-\alpha \ & \phantom{-}\beta \\
	\phantom{-}\alpha \ & -\beta
	\end{smallmatrix}
	\right]  ,  \left[
	\begin{smallmatrix}
	-\gamma \ & \phantom{-}\delta \\
	\phantom{-}\gamma \ & -\delta
	\end{smallmatrix}
	\right] \right\}$ for some $\alpha, \beta, \gamma, \delta > 0$. 
	This system is column weakly diagonally dominant, hence $\|x\|_1$ is a \gls{wPLF} \cite{BM15}.
	If we take $x(0)$ such that $x_1(0) + x_2(0) =1$, then $x_1(t) + x_2(t) = 1$ for all $t$.
However, on the plane  $x_1 + x_2 = 1$, $x_2(t)$ is not constant under persistent switching.
\hfill $\square$
\end{example}
}
 
\blue{
\begin{example}
For the \gls{CT} system with $A(w) =
	\left[
	\begin{smallmatrix}
	0 \ & -w \\
	w \ & \phantom{-}0
	\end{smallmatrix}
	\right],$ with $w > 0$, $V(x) = \|x\|^2$ is a \gls{wQLF}. While $\dot{V}(x) = 0$, the system exhibits persistent oscillations since the eigenvalues are on the imaginary axis.
\hfill $\square$
\end{example}
}

\blue{
Moreover, while (marginal) stability is equivalent to the existence of a homogeneous Lyapunov function \cite{chesi:03},  
Theorem \ref{th:strongconv-wplf_differentialinclusion} and the previous examples show that the existence of a wPLF is necessary but not sufficient for convergence. In the next statement, we show instead a necessary and sufficient condition, which is the existence of a wPLF with special conditions.}

\begin{corollary}
	\label{th:strongnsconv-wplf_differenceinclusion}
	The difference inclusion in \eqref{eq:polytopicdiffincl} (resp., differential inclusion in \eqref{eq:polytopicdiffinclcont}) is strongly convergent if and only if there exist a full row rank matrix $X$ and matrices 
	$$
	P_i=\left[\begin{array}{ccc}
	P_{i}^{\textup{as}} & & 0\\
	P_{i}^{\textup{r}} & & I_m
	\end{array}
	\right] \;
	\left(\text{resp., }
	P_i=\left[\begin{array}{ccc}
		P_{i}^{\textup{as}} & & 0\\
		P_{i}^{\textup{r}} & & 0_m
	\end{array}
	\right] \right),
	$$
	where $\{P_{i}^{\textup{as}}\}_{i \in \mc{M}}$ are such that $\|P_{i}^{\textup{as}}\|_1<1$ (resp., $\mu(P_{i}^{\textup{as}})<0$) and, for all $i \in \mc{M}$,
	$$
	A_i	 X=X P_i.
		\vspace{-0.75cm}
	$$
{\hfill $\square$}
\end{corollary}

	\vspace{-1.2cm}

\blue{
\begin{pf}
	We give the proof for the \gls{CT} case, since the one in \gls{DT} is analogous. The necessity part (only if) has been proved in Theorem~\ref{th:strongconv-wplf_differentialinclusion}. 
For the sufficiency (if), we consider the augmented system $\dot z=P(w)z$, with $z(0)$ such that $x_0=Xz(0)$, so that $x(t)=Xz(t), t\ge 0$. Such a system is strongly convergent in view of Theorem \ref{th:differentialinclusions_strongconv}:
$z(t) \rightarrow \bar z$, with $\bar z$ in the common kernel of the $P_i$, i.e. $P_i \bar z =0$.
We conclude the proof by noticing that $x$ converges as well: $x(t) = X  z(t) \rightarrow X \bar z \eqqcolon \bar x$, which belongs to $\textrm{ker}\left(\left\{A_i\right\}_{i \in \mc{M}}\right)$. The last claim is immediate, because $A_i \bar x = A_i X \bar z =   X P_i \bar z =0$.	{\hfill$\blacksquare$}
\end{pf}
}


\vspace{-0.25cm}

As for LTI systems, let us give (sufficient) LMI conditions to characterize strong convergence of difference and differential inclusions.   

\begin{theorem}\label{th:strongquadraticLMI}
	The difference inclusion in \eqref{eq:polytopicdiffincl} (resp., differential inclusion in \eqref{eq:polytopicdiffinclcont}) is strongly convergent if there exist $P \succcurlyeq 0$, with ${\rm rank}(P)=n-{\rm dim}(\bar {\mathcal X})$, and $Q\succ 0$ such that, for all $i \in \mc{M}$,
	\begin{equation}\label{eq.nuova}
	\begin{aligned}
	&\;\blue{A_i^\top P A_i - P + (A_i^\top-I)  Q  (A_i-I) \preccurlyeq 0}\\
	&\\
	&\left(\text{resp.,}\, \blue{A_i^\top P+PA_i +A_i^\top Q A_i \preccurlyeq 0}\right).
	\end{aligned}
	\vspace{-0.25cm}
	\end{equation}
{\hfill$\square$}
\end{theorem}
	\vspace{-0.5cm}

\begin{pf}
	\blue{ The proof in \gls{CT} mimics the one for the implication $(g) \Rightarrow (d)$ of Proposition~\ref{th:LTI_CT}, by replacing $A$ with $A_i$ and $\hat A$ with $\hat A_{i1}$, for all $i \in \mc{M}$.  The positive definite matrix $P_1$ generates a common \gls{wQLF} for the subsystem with matrices $\{\hat A_{1i}\}_{i \in \mc{M}}$ and hence strong convergence. 	Analogously, the proof in the DT case mimics the one for the implication $(d) \Rightarrow (f)$ of Proposition~\ref{th:LTI_DT}.  
}
	{\hfill$\blacksquare$}
\end{pf}

Finally, we characterize the relation between the strong convergence of the differential inclusion in \eqref{eq:polytopicdiffinclcont} and the \blue{associated} Euler difference inclusion 
\begin{equation}
\label{eq:nuovanuova}
x(k+1)=(I+ \tau A(w(k)) )x(k).
\end{equation}
\begin{proposition}
	\label{th:nuovonuovo}
	The differential inclusion in \eqref{eq:polytopicdiffinclcont} is strongly convergent if and only if there exists $\tau >0$ such that the  difference inclusion in \eqref{eq:nuovanuova} is strongly convergent. 
{\hfill$\square$}
\end{proposition}
	\vspace{-0.5cm}
\begin{pf}
	For all $i \in \mc{M}$, let $A_{\tau,i} \coloneqq I+\tau A_i$. The proof follows from Theorem~\ref{th:strongnsconv-wplf_differenceinclusion} by first noticing that ${\rm ker} (A_i) = {\rm ker} (A_{\tau,i}-I)$, for all $i \in \mc{M}$. Moreover, the fact that $\mu(P_{i}^{\textup{as}})<0$ implies $\|I+\tau P_{i}^{\textup{as}}\|_1<1$ for all $\tau <\min_{k} \, \{-2/d_{kk}\}$, where $d_{kk}$ are the diagonal entries of  $P_{i}^{\textup{as}}$. Finally, if $\tau >0$ is such that $\| P_{i}^{\textup{as}}\|_1<1$, then $\mu_1((P_{i}^{\textup{as}}-I)/\tau)<0$.
	{\hfill$\blacksquare$}
\end{pf}

\blue{
\subsection{Convergence rate}
We close the section with a convergence rate estimate for strong convergence. Let us consider the decomposition in Theorem~\ref{th:differentialinclusions_strongconv}. Since $A^\textup{as}(w)$ defines an exponentially stable inclusion, the first variable $x_1(t)$ converges to zero exponentially fast, i.e., $\| x_1(t) \| \leq c_0 \| x_1(0) \| e^{-\beta t}$, for some $0 < \beta := - \max_{i \in \mathcal{M}} \mu(P^\textup{as}_i)$, see Corollary \ref{th:strongnsconv-wplf_differenceinclusion}. The second component satisfies $\dot x_2 = A^\textup{r}(w) x_1$, thus 
	$$
	x_2(t) = x_2(0) + \int_0^t A^\textup{r}(w(\tau)) x_1(\tau) d \tau.
	$$
Consider $c_1$ such that $\| A^\textup{r}(w)\| \leq c_1$ for all $w$. Then, we have the following convergence rate for $x_2(t)$:
\begin{multline}
\| x_2(t) - x_2(\infty) \| = \|  \int_t^\infty  A^{ \textup{r} }(w(\tau)) x_1(\tau)  d \tau\| \\
		 \leq     
		c_1  \int_t^\infty  \|x_1(\tau)  \|d \tau  
		   \leq  
		c_0 c_1 \int_t^\infty  e^{-\beta t} d \tau 
		= \frac{c_0 c_1}{\beta}e^{-\beta t}.
\end{multline}
}

\section{Weak convergence}

\subsection{Lyapunov-like results}
We now characterize weak convergence \blue{for} uncertain systems \blue{via} Lyapunov arguments. \blue{First, we have a converse Lyapunov theorem.}

\begin{theorem}\label{th:differentinclusions_weakconv_cvx}
	If the difference inclusion in \eqref{eq:polytopicdiffincl} (resp. the diffrential inclusion \eqref{eq:polytopicdiffinclcont}) is weakly convergent, then it admits a \gls{wCLF}.
	{\hfill $\square$}
\end{theorem}
	\vspace{-0.5cm}
\begin{pf}
	Let us consider any polytopic, $0$-symmetric, set $\mc{X} \subseteq \R^n$ including the origin
	in its interior, and denote $\mc{R}_k$ as the set of reachable states for the uncertain system \eqref{eq:polytopicdiffincl} in at most $k$ steps from  $x(0) \in \mc{X}_0$.
	Let $\mc{R}_\infty$ be the union of all these sets, i.e., $\mc{R}_\infty \coloneqq \bigcup_{k \in \N}~\mc{R}_k$.
	By \cite{BS79}, $\mc{R}_\infty$ is bounded, namely there exists some $b \geq 0$ such that $\|x(k)\| \leq b$, for all $x \in \mc{R}_\infty$.
	Indeed, for all $k \in \N$, there exists some $x(k)$ that belongs to a trajectory
	originating from $x(0) \in \mc{X}_0$ (actually on a vertex \cite{BS79})
	which is on the boundary. Therefore, the boundedness of $\mc{R}_\infty$ is necessary for weak convergence.
	Thus, by defining $ \mc{C} \coloneqq \conv\{\mc{R}_\infty\}$, we note that its closure is a convex, compact, $0$-symmetric set including the origin in its interior, and by \cite[Th.~5.3]{B91} the norm induced by this set is a \gls{wCLF}.
	 The proof is similar (although a bit more involved) in the \gls{CT} case.
	{\hfill$\blacksquare$}
\end{pf}

	\vspace{-0.5cm}

\blue{Then, we provide sufficient \blue{matrix inequality} conditions for weak convergence via a common wQLF.}
  
\begin{theorem}\label{th:weakconv-cqlf_differentialinclusion}
	\blue{If there exist $P\succ 0$ and  $\eta\in (0,\, 1)$ (\text{resp., $\epsilon>0$ }) such that, for all $i \in \mc{M}$, the LMI
	\begin{equation}\label{eq:quadraticdt_LMI}
	\begin{aligned}
	&\; \eta(A_i^\top PA_i -P)+ (1-\eta) (A_i^\top - I) P (A_i -I) \preccurlyeq 0\\
	&\\
	&\left(\text{resp., }
		A_i^\top P+PA_i + \epsilon A_i^\top P A_i \preccurlyeq 0\right)
			\end{aligned}
	\end{equation}
}
	holds, then the difference inclusion in \eqref{eq:polytopicdiffincl} (resp., differential inclusion in \eqref{eq:polytopicdiffinclcont}) is weakly convergent.
{\hfill $\square$}
\end{theorem}
	\vspace{-0.5cm}
\begin{pf}
	Let $\epsilon=(1-\eta)/\eta$. By considering the \gls{wQLF} $V(x)=x^\top P x$, the Schur complement applied to the first \blue{inequality} in \eqref{eq:quadraticdt_LMI} gives, for any $w \in \mc{W}$,
	$$
	\left[\begin{array}{ccc}
	V(x(k+1))-V(x(k)) & & \vphantom{-} (x(k+1)^\top - x(k)^\top) P\\
	P (x(k+1)-x(k)) & & -\tfrac{1}{\epsilon}P
	\end{array}\right]\preccurlyeq 0.
	$$
	Then, $V(x(k+1))-V(x(k))\le 0$ and $x(k+1)=x(k)$  for some $x(k)= \bar{x}$. By noticing that the first \blue{inequality} in \eqref{eq:quadraticdt_LMI} can be equivalently rewritten as 
	$$
	\left[\begin{array}{ccc}
	P & A(w(k))^\top P& (A(w(k))^\top - I)P \\
	P A(w(k)) & P & 0\\
	P(A(w(k))- I) & 0 &  \tfrac{1}{\epsilon}P
	\end{array}\right]\succcurlyeq 0,
	$$
	for all $k \geq 0$, and by averaging in the interval $[0, K-1]$, we finally conclude that $\bar x$ belongs to $\ker(A(\bar w))$, where 
	$\textstyle
	\bar w=\lim_{K\to \infty}\frac{1}{K} \sum_{k = 0}^{K-1} w(k). 
	$
	Hence, the difference inclusion in \eqref{eq:polytopicdiffincl} converges and $V(x)$ is a common \gls{wQLF}.
	
	The \gls{CT} version is analogous: we consider the \gls{wQLF} $V(x)=x^\top P x$ and we apply the Schur complement to the second \blue{inequality} in \eqref{eq:quadraticdt_LMI}, hence obtaining, for any $w \in \mc{W}$,
	$$
	\left[\begin{array}{ccc}
	\dot V(x) & & \vphantom{-}\dot x^\top P\\
	P \dot x & & -\tfrac{1}{\epsilon}P
	\end{array}\right]\preccurlyeq 0.
	$$
	Then, $\dot V(x)\le 0$ and $\dot V(x)=0$ for $\dot x(t)=0$, i.e, for some $x = \bar{x}$. By noticing that the second \blue{inequality} in \eqref{eq:quadraticdt_LMI} can be equivalently rewritten as 
	$$
	\left[\begin{array}{ccc}
	A(w(t))^\top P+PA(w(t)) & & A(w(t))^\top P\\
	PA(w(t)) & & -\tfrac{1}{\epsilon}P
	\end{array}\right]\preccurlyeq 0,
	$$
	for all $t \geq 0$, and by averaging in the interval $[0, T]$, we finally conclude that $\bar x$ is in the null space of $A(\bar w)$, where 
	$\textstyle
	\bar w=\lim_{T\to \infty}\frac{1}{T} \int_0^T w(t) \; dt.
	$
	Thus, the differential inclusion converges and $V(x)$ is a common \gls{wQLF}.
	{\hfill$\blacksquare$}	
\end{pf}

\begin{remark}
The first \blue{inequality} in \eqref{eq:quadraticdt_LMI} (resp., the second \blue{inequality} in \eqref{eq:quadraticdt_LMI}) is equivalent to 
\begin{equation*}
\left(A_{\eta, i}^\textup{dt}\right)^\top P A_{\eta, i}^\textup{dt} \preccurlyeq P \quad \left(\text{resp.,} \left(A_{\epsilon,i}^\textup{ct}\right)^\top P + PA_{\epsilon,i}^\textup{ct} \preccurlyeq 0\right)
\end{equation*}
for all $i \in \mc{M}$, where $A_{\eta,i}^\textup{dt} = \eta A_i - \tfrac{1-\eta}{\eta} I$ (resp., $A_{\epsilon, i}^\textup{ct}=A_i(I+\epsilon)^{-1}$), 
implying weak stability of the auxiliary difference inclusion (resp., differential inclusion) 
$$
x(k+1) = A_{\eta}^\textup{dt} (w(k)) x(k), \; \left(\text{resp., }\dot x(t) = A_{\epsilon}^\textup{ct}(w(t)) x(t)\right).
$$
{\hfill$\square$}
\end{remark}

\blue{
\begin{remark}
	Differently from strong convergence, we cannot characterize the speed of convergence in the weak case, since it can be arbitrarily slow and strongly depends on the uncertain parameter. An example in \gls{CT} is the scalar system $\dot x = w(t) \, x$, with $w(t) \in [0, 1]$, for each $t \geq 0$.
{\hfill$\square$}
\end{remark}
}

\subsection{LaSalle-like criteria}
A necessary condition for convergence requires that the difference (differential) inclusion is at least marginally stable. To have convergence to a non-zero vector $\bar x$, at least one matrix inside the convex hull of the matrices $\{A_i-I\}_{i \in \mc{M}}$ (resp., $\{A_i\}_{i \in \mc{M}}$) shall be singular. In this case, the convergence does not need a common kernel and it can be associated with a specific singular element. As an example in \gls{CT}, let us consider
$$  
A(w(t)) \in \left\{A_1, A_2\right\} = \left\{ 
\left[ 
\begin{smallmatrix}
0 & & -1\\
1 & & -1
\end{smallmatrix}
\right], 
\left[ 
\begin{smallmatrix}
\phantom{-}0 & & \phantom{-}1 \\
-1& & -1
\end{smallmatrix}
\right]\right\},
$$
where $\bar A= \tfrac{1}{2}(A_1+A_2)$ is singular. Then, if $A(w)$ ``converges" to  $\bar A$, \blue{then} the trajectory converges to a non-zero vector of the form $[\bar{x}_1; 0]$, otherwise $x \to 0$. \blue{Next, we investigate LaSalle arguments of this type.}

\begin{lemma}\label{th:phi}
	Let $V$ be a positively homogeneous \gls{wLF} for the differential inclusion in \eqref{eq:polytopicdiffinclcont} with derivative
	$$
	D^+V(x,A(w)x) \coloneqq \underset{h\rightarrow 0^+}{\mathrm{lim \, sup}} \, \frac{V(x+hA(w)x)-V(x)}{h},
	$$
	and assume that
	\begin{equation}
	D^+V(x,A(w)x) \leq -\phi(x, A(w)x), \label{funphi}
	\end{equation}
	for some positive semi-definite, locally Lipschitz function $\phi$. Then $x(t)$, solution to \eqref{eq:polytopicdiffinclcont}, converges to the set
	\begin{equation*}
		{\mathcal N} \coloneqq \left \{ x \in \R^n \, \Big| \, \min_{w \in \mc{W}} \{\phi(x, A(w)x) \} = 0 \right \}.
	\end{equation*}
{\hfill$\square$}
\end{lemma}

\begin{pf} Let $x(t)$ be any trajectory of \eqref{eq:polytopicdiffinclcont}. We have:
	\begin{align}\label{eq:phi_proof}
	V(x(0)) & \geq  \int_0^t - D^+ V \, d\sigma \nonumber
\geq \int_0^t  \phi(x(\sigma), A(w(\sigma)) x(\sigma)) \, d\sigma \nonumber\\
	& \geq \int_0^t \min_{w \in \mc{W}} \, \{\phi(x(\sigma), A(w)x(\sigma)) \} \, d\sigma .
	\end{align}
	The function $\phi$ is Lipschitz continuous inside
	the set $\{x \in \R^n \mid V(x) \leq V(x(0)\}$. Since $\mathcal{W}$ is a compact and convex set, the min function is also nonnegative and Lipschitz \cite[Th.~7, p.~93]{Aub}. Therefore, in view of the Barbalat's lemma, ${\lim}_{t\to\infty} \, {\min}_{{w \in \mc{W}}} \, \{\phi(x, A(w)x) \} = 0$.
	{\hfill$\blacksquare$}
\end{pf}
	\vspace{-0.5cm}
	
The next two results follow as a direct consequence of Lemma~\ref{th:phi}.

\begin{proposition}\label{th:barba1}
	Let $V$ be a smooth \gls{wLF} for the differential inclusion in \eqref{eq:polytopicdiffinclcont}. Then $x(t)$, solution to \eqref{eq:polytopicdiffinclcont}, converges to the set
	\begin{equation}\label{eq:N_smooth}
	{\mathcal N} \coloneqq \left \{ x \in \R^n \, \Big| \, \min_{i \in \mc{M}} \, \{ - \nabla V(x) A_i x\} = 0 \right \}.
		\vspace{-0.5cm}
	\end{equation}
{\hfill$\square$}
\end{proposition}
	\vspace{-0.5cm}
	
\begin{pf}
	The proof follows by applying Lemma~\ref{th:phi} with $\phi(x, A(w)x)  = -\nabla V(x) A(w) x \geq 0$ and by noticing that the minimum in \eqref{eq:phi_proof} is achieved on the vertices, i.e.,
	$\underset{w \in \mc{W}}{\min} \, \{-\nabla V(x) A(w) x\}=  \underset{i \in \mc{M}}{\min} \, \{-\nabla V(x) A_i x\}$.
	{\hfill$\blacksquare$}
\end{pf}

\begin{remark}
Barbalat's lemma requires Lipschitz continuity of $D^+V(x,A(w)x)$. Consequently, if the differential inclusion in \eqref{eq:polytopicdiffinclcont} admits a non-smooth \gls{wLF}, $V$, then we cannot guarantee the convergence to the set where ${\min} \, \{-D^+ V(x, A(w) x)\} = 0.$ In that case, we shall rely on some locally Lipschitz function $\phi$, as in Lemma~\ref{th:phi}.
{\hfill$\square$}
\end{remark}

\begin{proposition} \label{cor:EAS}
	Let $V$ be a positively homogeneous \gls{wCLF} for the Euler auxiliary difference inclusion in \eqref{eq:nuovanuova} and define
	$$
	\Delta_\tau V(x,w) \coloneqq \frac{V(x+\tau A(w)x)-V(x)}{\tau},
	$$
	which is non-positive. Then $x(t)$, solution to \eqref{eq:polytopicdiffinclcont}, converges to the set 
	$$
	{\mathcal D} \coloneqq \left \{ x \in \R^n \, \Big| \, \min_{w \in \mc{W}}~\{- \Delta_\tau V(x,w) \} = 0 \right \}.
	\vspace{-0.75cm}
	$$
{\hfill$\square$}
\end{proposition}
	\vspace{-0.5cm}
\begin{pf}
\blue{
Since $V$ is convex, for all $x$ and $w$, we have that $D^+V(x,A(w)x) \leq  \Delta_\tau V(x,w)$ \cite{RouHabLal77}. The proof then follows by taking  $-\phi=\Delta_\tau V$, which is Lipschitz continuous.
}
{\hfill$\blacksquare$}
\end{pf}

\blue{
Let us observe that if $V$ is a positively homogeneous \gls{wPLF} for \eqref{eq:polytopicdiffinclcont}, we can choose $\phi(x,A(w)x) = \Delta_\tau V(x,w)$ in \eqref{funphi}, for some small enough $\tau > 0$. This follows by the fact that a polyhedral function is a \gls{wPLF} for the differential inclusion \eqref{eq:polytopicdiffinclcont} if and only if, for some $\tau >0$, it is a \gls{wPLF} for the associated Euler auxiliary difference inclusion.
}


In the DT case, we have the following LaSalle-like statement.

\begin{proposition}
	Let $V$ be a positively homogeneous \gls{wCLF} for the difference inclusion \eqref{eq:polytopicdiffincl}.
	Then $x(k)$, solution to \eqref{eq:polytopicdiffincl}, converges to the set 
	$$
	{\mathcal D} \coloneqq \left \{ x \in \R^n \, \Big| \, \min_{w \in \mc{W}} \{ V(A(w) x) - V(x) \} = 0 \right \}.
		\vspace{-0.75cm}
	$$
{\hfill$\square$}
\end{proposition}

\begin{pf}
	We have
		$\textstyle V(x(0)) =  \sum_{k=0}^\infty~ [ V(x(k)) - V(A(w(k)) x(k))]
		\geq  \sum_{k=0}^\infty~ \min_{w \in \mc{W}}\{-\Delta V(x(k),w)\}$, where $-\Delta V(x(k),w(k)) = V(x(k)) - V(A(w(k)) x(k))$.
	The terms of the series are positive, hence boundedness implies that
	$\underset{w \in \mc{W}}{\min} \, \{-\Delta V(x(k),w)\} \rightarrow 0$.
	{\hfill$\blacksquare$}
\end{pf}
	\vspace{-0.5cm}

\subsection{Weak kernel}

To further characterize weak convergence for differential inclusions, \blue{let us introduce the notion of weak kernel.}

\begin{definition}\textup{(\textbf{Weak Kernel})}\\
The weak kernel of the differential inclusion in \eqref{eq:polytopicdiffinclcont} is denoted by the set
$$
\mc{K} \coloneqq \left\{ x \in \R^n \mid \exists w \in \mc{W} \; \text{ s.t. } \; A(w) x = 0 \right\}.
	\vspace{-0.75cm}
$$
\hfill$\square$
\end{definition}

By referring to the set $\mc{N}$ in \eqref{eq:N_smooth}, in general $\mc{K} \subseteq {\mathcal N}$. The next example shows a case in which $\mc{K}$ and $\mc{N}$ coincide, and the solution to \eqref{eq:polytopicdiffinclcont} converges to the weak kernel.

\begin{example}\label{exa:ker}
	For the differential inclusion with
	$	
	A(w(t)) \in \left\{A_1, A_2\right\} = \left\{ 
	\mathrm{diag}(-1,0), \, 	\mathrm{diag}(0,-1) \right\},
	$
	$V(x) = \tfrac{1}{2}(x_1^2+ x_2^2)$ is a \gls{wQLF}. The set $\mc{N}$ in \eqref{eq:N_smooth} is
	characterized by the equation
	$$
	\min_{i \in \mc{M}} \, \{x^\top A_i x\} = \min \, \{x_1^2,x_2^2\} =0,
	$$
	which represents the (nonconvex) union of the two axes. Note that, in general, if there exists a \gls{wQLF} for the differential inclusion in \eqref{eq:polytopicdiffinclcont}, since $A_i P + PA_i = -Q_i$, with $Q_i \succcurlyeq 0$ for all $i \in \mc{M}$, the set $
	{\mathcal N} = \left\{ x \in \R^n \mid \min_{i \in \mc{M}} \, \{ x^\top Q_i x \} = 0 \right\}
	$ 
	can be computed as the union of the kernels.
	{\hfill$\square$}
\end{example}

\blue{
Our result is that, under the existence of a \gls{wPLF}, asymptotic stability of the differential inclusion in \eqref{eq:polytopicdiffinclcont} is equivalent to the fact that the weak kernel is trivial.
}

\begin{theorem}
	Let $V$ be a \gls{wPLF} for the differential inclusion in \eqref{eq:polytopicdiffinclcont}. Then $A(w)$ is robustly 
	non-singular, i.e., ${\mathcal K} = \{0\}$, if and only if the solution to \eqref{eq:polytopicdiffinclcont},
	$x(t)$, tends to $0$ for all $x(0)$ and $(w(t))_{t\geq0}$.
{\hfill$\square$}
\end{theorem}
	\vspace{-0.5cm}
\begin{pf}
	If $A(w)$ is singular for some $w \in \mc{W}$, then we can not have robust stability. Thus, we only need to prove the converse statement
	Let $A(w)$ be robustly non-singular and let $V$ be defined by $s$-planes, indexed in $\mc{S} \coloneqq \{1, \ldots, s\}$. Given some $\rho > 0$, we introduce the set
	$
	\mc{L}_{(V/\rho)} \coloneqq \{x \in \R^n \mid V(x) \leq \rho\} = \{x \in \R^n \mid f^\top_i x \leq \rho, \; \forall i \in \mc{S}\}
	$.
	Moreover, let $\mc{A}(x)$ be the set indexing the active planes, i.e., $\mc{A}(x) \coloneqq \{i \in \mc{S} \mid f^\top_i x = V(x)\}$ and let $\mathcal{F}_\mathcal{A}$ be any arbitrary $\ell$-dimensional face 
	$$
	\mathcal{F}_\mathcal{A} \coloneqq \{x \in \mc{L}_{(V/\rho)} \mid f^\top_i x =\rho, \; \forall i \in  \mathcal{E}\},
	$$
	where $\mc{E}$ indexes the ``active'' planes at $\mathrm{bdry}\left(\mc{L}_{(V/\rho)}\right)$. The remainder of the proof exploits \cite[Lemma~3]{BG18}. Thus, for some $\ell$-dimensional face, and hence for some indices $i \in \mathcal{A}(x)$, we shall have, for all $x \in \mathcal{F}_\mathcal{A}$ and $w \in \mc{W}$, $D^+V(x,A(w)x) = f^\top_i A(w) x < 0$.
	Otherwise, the existence of some $\bar{w} \in \mc{W}$ such that $f^\top_i A(\bar w) x = 0$ for all $i \in \mathcal{A}(x)$ implies the positive invariance of the \gls{LTI} system $\dot x = A(\bar w) x$, and hence the singularity of  $A(\bar w)$.
	Thus, let $x(t)$ be the solution to \eqref{eq:polytopicdiffinclcont} with initial condition $x(t_0)$. If $x(t_0) \in \mathcal{F}_\mathcal{A}$, by continuity there exists a right neighborhood $[t_0, \tau_1]$ of $t_0$, with $\tau_1 > 0$, such that the active set can not grow, i.e, $\mathcal{A}(x(t)) \subset \mathcal{A}(x(t_0))$ for all $t \in [t_0, \tau_1]$.
	By repeating the same argument, $x(t)$ shall reach a face of dimension $1$ where there is only a single active constraint, i.e., $\mc{A}(x(t)) = \{i\}$. Here, in view of the considerations above, $f^\top_i A(w) x <0$ for all $x \in \mc{F}_\mc{A}$ and $w \in \mc{W}$.
	Therefore, for any initial condition $x(t_0)$ such that $V(x(t_0)) = \rho$, in an arbitrary small neighborhood $[t_0, \tau]$ of $t_0$, for some $\tau >0$, we have $V(x(\tau)) < \rho$.
	Then, $V(x(t)) \geq 0$ is monotonically non-increasing along any $x(t)$, and hence it has a limit $V(x(t)) \to \bar{\rho}$ from above. From now on, the proof replicates the one in \cite[Proof of Th.~2]{BG18}, and hence is here omitted.
	%
	{\hfill$\blacksquare$}
\end{pf}

\blue{
Finally, we note that robust non-singularity, i.e., trivial weak kernel, along with the existence of a \gls{wLF} (which is not polyhedral), does not imply convergence. For example, consider the system
	$\dot{x} = \left [ \begin{smallmatrix} 0 & -1 \\
	1 &  ~0
	\end{smallmatrix}
	\right ] x$ and $V = \|x\|$.
}




\section{Convergence lacks of duality}

It is known that the stability (asymptotic or marginal) of $x(k+1) = A(w(k)) x(k)$  (resp., $\dot x(t) = A(w(t)) x(t)$) implies the stability of the dual system
$x(k+1) = A^\top(w(k)) x(k)$ (resp., $\dot x(t) = A^\top(w(t)) x(t)$). \blue{Similarly, it can be shown that duality holds in the LTI case, i.e., if $x(k+1) = A x(k)$  (resp., $\dot x(t) = A x(t)$) converges, then $x(k+1) = A^\top x(k)$  (resp., $\dot x(t) = A^\top x(t)$) converges as well. Interestingly, the convergence of an uncertain \blue{linear} system lacks of duality. }

\begin{proposition}
\blue{
The convergence of $x(k+1) = A(w(k)) x(k)$  (resp., $\dot x(t) = A(w(t)) x(t)$) does not imply the convergence of the dual system
$x(k+1) = A^\top(w(k)) x(k)$ (resp., $\dot x(t) = A^\top(w(t)) x(t)$).
}
{\hfill$\square$}
\end{proposition}
	\vspace{-0.5cm}
\begin{pf}
	The proof goes by means of two examples. We consider the \gls{DT} case first. 
	With $a \in (0,1)$, the system
	$$
	A(w(k)) \in \{A_1, A_2\} = \left\{\left[
	\begin{array}{ccc}
	a & & 1\\
	0 & & 1
	\end{array}
	\right],
	\left[
	\begin{array}{ccc}
	a & & 2\\
	0 & & 1
	\end{array}
	\right]\right\}
	$$
	is neither strong nor weak convergent. In fact, while $x_2(k) = x_2(0)$, for all $k > 0$, by taking alternatively $A(w(k)) = A_1$ for $k$ even and $A(w(k)) = A_2$ for $k$ odd, $x_1(k)$ exhibits persistent oscillations. However, the dual system is strongly convergent, since it is in the block-triangular form in Theorem~\ref{th:differentialinclusions_strongconv} $(b)$.
	
	In \gls{CT}, the differential inclusion characterized by
	$$ A(w(t)) \in \{A_1, A_2\} = \left\{ 
	\left[
	\begin{array}{ccc}
	-\alpha & & \phantom{-}\beta \\
	\phantom{-}\alpha & & -\beta
	\end{array}
	\right]  ,  \left[
	\begin{array}{ccc}
	-\gamma & & \phantom{-}\delta \\
	\phantom{-}\gamma & & -\delta
	\end{array}
	\right] \right\},$$
	for some $\alpha, \beta, \gamma, \delta > 0$, does not converge (see Example~\ref{ex:QLFPLF_not_conv}) while the dual system with $A^\top(w(t))$ is strongly convergent to $\kappa \bsone$, $\kappa \in \R$. It follows by considering $V(x) = |x_1 - x_2|$ as a common Lyapunov function.
	\hfill$\blacksquare$
\end{pf}
	\vspace{-0.5cm}

However, an intriguing fact is that weak quadratic Lyapunov functions provide some kind of duality.

\begin{theorem}
	Assume that the difference inclusion in \eqref{eq:polytopicdiffincl} (differential inclusion in \eqref{eq:polytopicdiffinclcont}) admits a \gls{wQLF}. Then, it is strongly convergent if and only if the dual system is strongly convergent.
	{\hfill $\square$}
\end{theorem}
	\vspace{-0.5cm}
\begin{pf}
	We prove the \gls{CT} version first. By applying the state decomposition of Theorem~\ref{th:differentialinclusions_strongconv}, we consider a \gls{wQLF} for \eqref{eq:polytopicdiffinclcont} such that, for all $i \in \mc{M}$,
	\begin{equation*}
	\begin{aligned}
	&\left[\begin{array}{cc}
	P & R\\
	R^\top & Q
	\end{array}\right] 
	\left[\begin{array}{ccc}
	A_i^{\textup{as}} & & 0\\
	A_i^{\textup{r}} & & 0_m
	\end{array}\right]
	+  
	\left[\begin{array}{ccc}
	(A_i^{\textup{as}})^\top & & (A_i^{\textup{r}})^\top\\
	0 & & 0_m
	\end{array}\right] \left[\begin{array}{cc}
	P & R\\
	R^\top & Q
	\end{array}\right]  =\\ 
	&\left[\begin{array}{ccc}
	P A_i^{\textup{as}} +(A_i^{\textup{as}})^\top P + RA_i^{\textup{r}} + (R A_i^{\textup{r}})^\top & & (*) \\
	R^\top A_i^{\textup{as}}   +  QA_i^{\textup{r}}   & &  0_m
	\end{array}\right]   \preccurlyeq 0,
	\end{aligned}
	\end{equation*}
	where $(*) = (A_i^{\textup{as}})^\top R   +  (A_i^{\textup{r}})^\top Q$. The relation above implies  
		$\textstyle
		R^\top A_i^{\textup{as}}   +  Q A_i^{\textup{r}} = 0,
		$
		for all $i \in \mc{M}$. Then, we use the matrices $R^\top$ and $Q$ to introduce an additional state transformation as follows
		$$
		\left[\begin{array}{cc}
		I & 0\\
		R^\top & Q
		\end{array}\right]
		\left[\begin{array}{ccc}
		A_i^{\textup{as}} & & 0\\
		A_i^{\textup{r}} & & 0_m
		\end{array}\right]\left[\begin{array}{cc}
		I & 0\\
		R^\top & Q
		\end{array}\right]^{-1} =
		\left[\begin{array}{ccc}
		A_i^{\textup{as}} & & 0\\
		0 & & 0_m
		\end{array}\right].
		$$
		Since $\left\{A_i^{\textup{as}}\right\}_{i \in \mc{M}}$ generate asymptotically stable differential inclusions, the same holds for $\left\{(A_i^{\textup{as}})^\top\right\}_{i \in \mc{M}}$. Hence the strong convergence of $\mathrm{diag}\left( A_i^{\textup{as}}, 0_m \right)$
		implies the strong convergence of the linear differential inclusion defined by $\left\{A_i^\top\right\}_{i \in \mc{M}}$, where
		$
		A_i^\top=\left[\begin{array}{ccc}
		(A_i^{\textup{as}})^\top  & & (A_i^{\textup{r}})^\top\\
		0 & & 0_m
		\end{array}\right].
		$
	The proof is similar for the difference inclusion in \eqref{eq:polytopicdiffincl}. Indeed, by considering a \gls{wQLF}, we have
	\begin{equation*}
	\left[\begin{array}{cc}
	(A_i^{\textup{as}})^\top &  (A_i^{\textup{r}})^\top\\
	0 & I_m
	\end{array}\right]\left[\begin{array}{cc}
	P & R\\
	R^\top & Q
	\end{array}\right]\, \left[\begin{array}{cc}
	A_i^{\textup{as}} & 0\\
	A_i^{\textup{r}} & I_m
	\end{array}\right] - \left[\begin{array}{cc}
	P & R\\
	R^\top & Q
	\end{array}\right]\,
	\preccurlyeq 0, 
	\end{equation*}
	which implies, for all $i \in \mc{M}$,
	$\textstyle
	R^\top A_i^{\textup{as}}   +  QA_i^{\textup{r}} - R^\top = 0.
	$
	Moreover, by adopting the same additional state transformation introduced above, we have
	$$
	\left[\begin{array}{cc}
	I & 0\\
	R^\top & Q
	\end{array}\right]
	\left[\begin{array}{ccc}
	A_i^{\textup{as}} & & 0\\
	A_i^{\textup{r}} & & I
	\end{array}\right]\left[\begin{array}{cc}
	I & 0\\
	R^\top & Q
	\end{array}\right]^{-1} =
	\left[\begin{array}{ccc}
	A_i^{\textup{as}} & & 0\\
	0 & & I_m
	\end{array}\right],
	$$
	where $\{A_i^{\textup{as}}\}_{i \in \mc{M}}$ generate an asymptotically stable difference inclusion. So is $\left\{(A_i^{\textup{as}})^\top\right\}_{i \in \mc{M}}$, and by following the same reasoning for the \gls{CT} case, we obtain the strong convergence of the dual system.
	\hfill$\blacksquare$
\end{pf}

	\vspace{-0.5cm}

\blue{
Let us summarize the main implications for linear differential/difference inclusions next, in Table \ref{tavola}.
\begin{table}[!htb]
$$
\begin{array}{ccccccc}
\text{AS} & \iff & \exists \, \text{PLF} & \iff & \exists \, \text{PLF\textup{*}} & \iff & \text{AS\textup{*}} \\
\Downarrow &   &   &   &  &   &  \Downarrow \\
\text{sCON} & \Longrightarrow & \exists \, \text{wPLF} & \iff &\exists \,  \text{wPLF\textup{*}} & \Longleftarrow & \text{sCON\textup{*}} \\
\Downarrow &   &   &   &  &   &  \Downarrow \\
\text{wCON} &  &   &  &  &   & \text{wCON\textup{*}}\\
\Downarrow &   &   &   &  &   &  \Downarrow \\
\text{MS} & \iff   & \exists \, \text{wCVXLF}  & \iff  & \exists \, \text{wCVXLF\textup{*}} & \iff  & \text{MS\textup{*}}
\end{array}
$$
\caption{
Summary of the implications. AS = asymptotic stability;
PLF = (existence of) polyhedral Lyapunov function; sCON = strong convergence;  wCON = weak convergence; MS = marginal stability; wCVXLF = weak convex Lyapunov function. The asterisk $^*$ refers to the dual system.
}
\label{tavola}
\end{table}
}

\section{Examples}
\subsection{Uncertain linear systems with persistent input}
Let us consider the system
$$
\dot x(t) = A(w(t)) x(t) + B\bar u 
$$
where $\bar u \in \R$ is a constant scalar, and assume that the unforced system is asymptotically stable for any signal $w(t) \in \mc{W}$.
We observe that the state of such system converges if and only if the following extended (marginally stable) system is strongly convergent:
\begin{equation}
\label{eq:extended_system}
	 \left[\begin{array}{c}
	\dot x \\
	\dot u
	\end{array}\right] =  \left[\begin{array}{cc}
	A(w) & B\\
	\bs{0}_n^\top & 0
	\end{array}\right]
 \left[\begin{array}{c}
	 x \\
	 u
	\end{array}\right].
\end{equation}
In fact, strong convergence implies that the extended state converges to the common kernel $\mc{K} = \left\{ [\bar x \, ; \, \bar u] \in \R^{n+1} \mid A(w) \bar x + B \bar u =0 \right\}$. 
This subspace has dimension $1$, since $A(w)$ is non-singular. 
Let $( \bar{x}, \bar{u} ) \in \mathcal{K}$: with the transformation 
$
T=\left[ \begin{array}{ccc} I_{n} & & \bar x \\ \bs{0}_n^\top & & \bar u \end{array}\right],
$
we have that
$
	 T^{-1} \left[\begin{array}{cc}
	A(w) & B\\
	\bs{0}_n^\top & 0
	\end{array}\right]T= \mathrm{diag}\left(A(w),0 \right)	,
$
hence, by Theorem \ref{th:differentialinclusions_strongconv}, the extended system in \eqref{eq:extended_system} is strongly convergent. 

The DT case is analogous, since the uncertain affine system $x(k+1) = A(w(k)) x(k) + B \bar{u}$, where $\bar{u} \in \R$ is constant, reads as the following uncertain linear system, where $u(0) = \bar{u}$:
$$
\left[
\begin{matrix}
x^+ \\
u^+
\end{matrix}
\right] = \left[ 
\begin{matrix}
A(w) & B \\
\bs{0}_n^\top & 1
\end{matrix}
\right] \left[
\begin{matrix}
x \\
u
\end{matrix}
\right].
$$

\subsection{Uncertain opinion dynamics in social networks}

An example of uncertain system with persistent input is the following model of the probability distribution of a linear emulative network of stochastic agents. Let $x^{[r]}(t)\in \R^{M}$ denote the unit-sum opinion probability vector for the $r$-th agent, $r \in \mc{N}\coloneqq\{1,2,\ldots,N\}$, in linear emulative model with an arbitrary network topology. It is shown in \cite{Bolzern2019} that the probability distribution of the $r$-th agent can be written as, for $j \in \mc{M} \coloneqq \{1,2,\ldots,M\}$,
$$
\dot x_j^{[r]}=\sum_{i \in \mc{M}} q_{i,j}^{[r]}x_i^{[r]}+\lambda \left( \frac{\sum_{k\in {\mathcal N}_r} x_j^{[k]}}{{|{\mathcal N}_r|}}-x_j^{[r]}\right)
$$
where
$\lambda>0$ represents the influence intensity, equal for different opinions, ${\mathcal N}_r$ the set of neighbors of agent $r$, and  $q_{i,j}^{[r]}$  the probability transition rate from opinion $i$ to opinion $j$ for agent $r$ in absence of influence. From now on, we assume that all the transition rate matrices $Q^{[r]}$, $r \in \mc{N}$, of the isolated agents are irreducible and $M=2$. Since ${\bf 1}_2^\top x^{[r]}=1$, we can discard the second component of each probability vector, so that the reduced $N$-dimensional model reads as
$$
\dot z(t)=(-\Delta + \lambda G)z(t)+ q
$$
where $G$ is the graph matrix, a Metzler matrix such that ${\bf 1}_N^\top G = \bs{0}_N^\top$, $\Delta = \mathrm{diag}(\delta_1, \ldots, \delta_N)$ is a positive diagonal matrix with elements $\delta_r = q_{1,2}^{[r]}+q_{2,1}^{[r]}$ and $q$ is a column vector with $q_r=q_{2,1}^{[r]}$. Let us consider that the graph matrix is uncertain and given by
$$
\textstyle 
G(w)=\sum_{i=1}^g G^{[i]} w_i
$$
with $G^{[i]}$, $i=1,2,\ldots, g$, suitable graph matrices of connected graphs. 
The augmented system is then
\begin{equation} \label{eq:augmented_opinion_dyn}
\left[ \begin{array}{c} \dot z\\ \dot u\end{array}\right]=\left[ \begin{array}{ccc}\ \sum_{i=1}^g (-\Delta + \lambda G^{[i]})w_i & & q \\ 0 & & 0\end{array}\right] \left[ \begin{array}{c} z \\ u\end{array}\right].
\end{equation}
We note that no common kernel exists unless $q=\bar \beta\Delta {\bf 1}_N$, and in this case, the common ($1$-dimensional) kernel is given by 
$\mathrm{span}( [ {\bf{1}}_N^\top, \,  ( \Delta \, {\bf{1}}_N)^\top ]^\top )$,
since 
$-(-\Delta + \lambda G^{[i]})^{-1}{\bf 1}_N={\bf 1}_N$, $\forall i \in \mathcal{N}$.
This is the well known {\it unbiased} network case, i.e., the bias parameters $\beta_r =
{q_{21}^{[r]}}/({q_{21}^{[r]}+q_{12}^{[r]}})$ are all equal to $\bar\beta<1$. 
It follows from \S 8.1 that the augmented uncertain system in \eqref{eq:augmented_opinion_dyn} is strongly convergent, in particular, for all $z(0)$ and $u(0)$, $\lim_{t \rightarrow \infty }z(t) = {\bf 1}_N \bar\beta u(0)$.

\subsection{Kolmogorov--like equations}

Let us consider the system 
\begin{equation}\label{eq:kolm_like}
\dot x(t) = \left( \textstyle \sum_{i \in \mc{M}}  A_i w_i(t) \right) x(t),
\end{equation}
where each $A_i$ is an irreducible, Metzler matrix with strictly positive off-diagonal elements. 
If $A_i$ are column-stochastic matrices, i.e., if $\bs{1}^\top A_i = \bs{1}^\top$ for all $i \in \mc{M}$, then $V(x) = \|x\|_{1}$ is a \gls{wPLF} for the system in \eqref{eq:kolm_like}, while if $A_i$ are row-stochastic matrices, i.e., if $A_i \bs{1} = \bs{1}$ for all $i \in \mc{M}$, then $V(x) = \|x\|_{\infty}$ is a \gls{wPLF}.

In general, as shown in Example~\ref{ex:QLFPLF_not_conv}, column-stochastic matrices do not share a common kernel. Conversely, in the row-stochastic case, for $\tau > 0$ small enough, we note that the associated Euler difference inclusion
$
x(k+1) = F(w(k)) x(k)= \left(\sum_{i \in \mc{M}} w_i (k) F_i \right) x(k),
$
where  $F_i  \coloneqq I+\tau A_i$, $\forall i \in \mc{M}$, is a positive, row-stochastic matrix. Then, we have
$
\| F_i \, x\|_\infty \leq \|x\|_\infty
$
where, unless $x_1=x_2=\dots =x_n$, the relation holds as strict inequality since, for all $i \in \mc{M}$, $\sum_j \left(F_i\right)_{hj} |x_j| < \max_j |x_j|$. Therefore, the state of the system converges to the set $\mc{K} =\left\{x \in \R^n \mid x_1=x_2=\dots =x_n\right\}$. Under the irreducibility assumption, $\mc{K}$ corresponds to the common kernel of the generator matrices $\{A_i\}_{i \in \mc{M}}$.

\subsection{Plant tuning with singularity}
Given a static plant governed by an \emph{unknown} mapping
$$
y(t) = f(u(t))
$$
such that the Jacobian matrix of $f$ belongs to a polytope  $\mc{M}$ in which all the elements are non-singular (robust non-singularity) and such that $f(\bar u) =0$ for some unique $\bar u$, then there exists a dynamic control law $\dot u(t) = \phi(y(t))$
that steers $u(t)$ to $\bar u$ \cite{blanchini2017model}.  Conversely, let us consider the case in which there are isolated singular elements in $\mc{M}$, hence we cannot guarantee convergence to $0$, but we still need to ensure state convergence.

As an example, let us consider the flow system governed  by the equations


\vspace{-0.5cm}

$$
\left\{
\begin{aligned}
&y_1 = -\phi(u_1 - u_2) +r_1\\
&y_2 = \phi(u_1 - u_2) + \psi(u_2) -r_1 +r_2,
\end{aligned}
\right.
$$

where $r_1$ and $r_2$ are references, $a_{\min} \leq a \coloneqq  \phi^{\prime}(u) \leq a_{\max}$, $b_{\min} \leq b \coloneqq \psi^{\prime}(u) \leq b_{\max}$. 
Next, we consider an integral control law $\dot{u}  = - k \,  y$, with $k > 0$ and derive the dynamics $\dot{y} = \mathrm{J}f(u) \, \dot{u}$, i.e.,
$$
 \left[\begin{array}{cc}
	~\dot y_1~ \\
	~\dot y_2~
	\end{array}\right]=   - k \left[\begin{array}{cc}
	\phantom{-}a  ~~ &  -a~ \\
	-a~ ~& ~~a+b~
	\end{array}\right]~
\left[\begin{array}{cc}
	~y_1~ \\
	~y_2~
	\end{array}\right]\,,
$$
\vspace{-0.5cm}

which is a linear differential inclusion since $a \in [a_{\min}, a_{\max}]$ and $b \in [b_{\min}, b_{\max}]$.
For this system, $V(y) = \| y \|$ is a \gls{wLF} and, in view of Corollary~\ref{th:barba1}, the state weakly converges to the set $\mc{N} = \left\{y \in \R^2 \mid y_1 - y_2 =0\right\}$.

\section{Conclusion}
The state convergence problem is highly relevant within the system-and-control community, since it occurs in several areas, from multi-agent learning, consensus and opinion dynamics to plant tuning. For \gls{LTI} systems, state convergence can be characterized via necessary and sufficient \blue{linear} matrix inequalities and Lyapunov-like conditions. 
In the presence of uncertainty, two different definitions of state convergence, i.e., strong and weak convergence, shall be considered. These two are equivalent under the kernel sharing property.
In general, while strong convergence is structurally guaranteed by the separability of the eigenvalues of the generator matrices, weak convergence is not. 
\blue{Lyapunov-like sufficient conditions for weak convergence can be established via linear matrix inequalities.} Necessary and sufficient conditions for weak convergence \blue{are currently unknown}.


\appendix

\section{Proofs of statements in Section~\ref{sec:LTI}}\label{sec:app_linear}

\subsection{Proof of Lemma \ref{lemma:convergent-stable}}
Let $\D_r \coloneqq \{ z \in \mathbb{C} \mid | z - (1-r) | \leq r \}$, the disk of radius $r>0$ centered in $(1-r,0)$. 
\blue{
The system in \eqref{eq:LTI_DT_aux} is (marginally) stable if and only if the eigenvalues of $A^{\textup{dt}}_\eta$ are contained in the unit disk $\D_1$ and those at $1$ are semi-simple. By \cite[Lemma~4]{BFBG18} the eigenvalues of $A$ are in $\D_\eta$ and those at $\mathrm{bdry}(\D_\eta)$, which always containts $\{ 1\}$, are semi-simple. Therefore, the system in \eqref{eq:LTI_DT} is (marginally) stable and it converges.
}
\hfill$\blacksquare$

\subsection{Proof of Proposition \ref{th:LTI_DT}}
	$(a) \Leftrightarrow (b)$ follows from Lemma~\ref{lemma:convergent-stable}, while $(b) \Leftrightarrow (c)$ follows by \cite[Th.~4.50]{BM15}.
	$(a) \Leftrightarrow (d)$ \blue{follows} from the fact that $A$ has all the eigenvalues strictly inside the unit disk and the eigenvalues in $1$ are semi-simple, \blue{thus} there always exists a state transformation that allows to separate stable and critical eigenvalues. 
	$(c) \Leftrightarrow (e)$ : By the Lyapunov theorem, the system in \eqref{eq:LTI_DT_aux} is \blue{(marginally)} stable if and only if
$
	\left({A^{\textup{dt}}_\eta}\right)^\top P A^{\textup{dt}}_\eta \preccurlyeq P$,
	\blue{which leads to inequalities (e). \\
	$(d) \Rightarrow (f):$ Assume that 
$$T^{-1} A T = 
		\left[ \begin{smallmatrix}
		A^{\textup{as}} & & 0 \\ 
		A^{\textup{r}} & & I 
		\end{smallmatrix} \right],
		$$
		for some Schur matrix $A^{\textup{as}}$;
and take $P_1\succ 0$ such that 
	$$
	(A^{\textup{as}})^\top P_1 A^{\textup{as}} -P_1+W\prec 0,
	$$
	where 
	$
	W \coloneqq (A^{\textup{as}}-I)^\top (A^{\textup{as}}-I) + (A^{\textup{r}})^\top A^{\textup{r}} \succ 0$. Then the inequality in the statement is satisfied with $P$ and $Q$ such that $T^\top PT={\rm diag }(P_1, \, 0)$ and $T^\top QT=I$. 
	\\
	$(f) \Rightarrow (d)$ :
	\sloppy Let $\hat A=U^\top Q^{1/2}AQ^{-1/2}U$, $\hat P=U^\top Q^{-1/2}PQ^{-1/2}U$, where $U$ is an orthogonal matrix such that         
	$
  	\hat P= \mathrm{diag}( P_1,0 ),
 	$
	with $P_1\succ 0$. Hence, the \blue{matrix inequality} in the statement can be rewritten as  
	$$
	\hat A^\top \hat P \hat A-\hat P+ (\hat A^\top -I)(\hat A-I) \preccurlyeq 0 \,,
	$$        
	Letting	$
	\hat  A= \left [
	\begin{smallmatrix} 
	\hat A_1 & &\hat A_3 \\ 
	\hat A_2 & & \hat A_4 
	\end{smallmatrix} \right ],
	$
the previous inequality implies that $\hat A_3=0$, $\hat A_4=I$, so that 
	$
	\hat  A= \left [
	\begin{smallmatrix} 
	\hat A_1 & & 0 \\ 
	\hat A_2 & & I 
	\end{smallmatrix} \right ],
	$. 
Moreover,  
	$$
	\blue{P_1 \succcurlyeq  \hat A_1^\top P_1 \hat A_1+  W},
	$$
where  
$
\hat  W \coloneqq (\hat A_1^\top -I)(\hat A_1-I) + \hat A_2^\top \hat A_2.
$
Due to the assumption on the rank of $P$, $\hat W$ is positive definite, hence $A_1$ is Schur and the system is convergent.
}
{\hfill$\blacksquare$}

\subsection{Proof of Lemma \ref{lemma:convergent-stable_tc}}

	The proof follows by noticing that there exists a relation among the eigenvalues $\lambda \in \Lambda(A)$ and $\nu \in \Lambda(A^{\textup{ct}}_\epsilon)$ that does not alter both geometric and algebraic multiplicities, i.e.,
	$
	\nu = (\lambda)/({1 + \epsilon \lambda}).
	$
	The system in \eqref{eq:LTI_CT} converges if and only if all its eigenvalues have strictly negative real part and those on the imaginary axis are semi-simple and all equal to $0$. Hence, there always exists $\epsilon > 0$ small enough that maps each $\lambda \in \Lambda(A)$ with strictly negative real part to some $\nu \in \Lambda(A^{\textup{ct}}_\epsilon)$ belonging to the open left-half plane. Therefore, by noticing that the eigenvalues of $A$ at $0$ are mapped in $0$, the system in \eqref{eq:LTI_CT_aux} is stable.
	\hfill$\blacksquare$

\subsection{Proof of Proposition \ref{th:LTI_CT}}

	$(a) \Leftrightarrow (b)$ follows from Lemma~\ref{lemma:convergent-stable_tc}.
	$(a) \Leftrightarrow (c) \Leftrightarrow (d)$ \blue{follows} from \cite[Th.~4.49]{BM15}.
	$(a) \Leftrightarrow (e)$ : As in Lemma~\ref{lemma:convergent-stable_tc}, there exists a relation among the eigenvalues $\lambda \in \Lambda(A)$ and $\nu \in \Lambda(A_\tau)$ that does not alter both geometric and algebraic multiplicities, i.e., $\nu = 1 + \tau \lambda$. Moreover, the eigenvectors of $A_\tau$ are respectively the same of those of $A$. Therefore, the eigenvalues at $0$ of $A$ are semi-simple if and only if the eigenvalues at $1$ of $A_\tau$ are semi-simple.\\
	$(b) \Leftrightarrow (f)$ : By the Lyapunov theorem, the system in \eqref{eq:LTI_CT_aux} is \blue{(marginally)} stable if and only if
	$
	\left({A^{\textup{ct}}_\epsilon}\right)^\top P + P A^{\textup{ct}}_\epsilon \preccurlyeq 0.
	 $
	Hence, by pre and post-multiplying by $(I + \epsilon A)^\top$ and $(I + \epsilon A)$, respectively, we obtain the \blue{inequality} in $(f)$.
\blue{
$(d) \Rightarrow (g):$ Assume that 
$$T^{-1} A T = 
		\left[ \begin{smallmatrix}
		A^{\textup{as}} & & 0 \\ 
		A^{\textup{r}} & & 0 
		\end{smallmatrix} \right],
		$$
		for some Hurwitz matrix $A^{\textup{as}}$;
and take $P_1\succ 0$ such that 
	$$
	(A^{\textup{as}})^\top P_1+ P_1 A^{\textup{as}} +W\prec 0,
	$$
	where 
	$
	W \coloneqq (A^{\textup{as}})^\top A^{\textup{as}} +( A^{\textup{r}})^\top A^{\textup{r}} \succ 0$. Then the inequality in the statement is satisfied with $P$ and $Q$ such that $T^\top PT={\rm diag }(P_1, \, 0)$ and $T^\top QT=I$. 
	\\
	$(g) \Rightarrow (d)$ :
	\sloppy Let $\hat A=U^\top Q^{1/2}AQ^{-1/2}U$, $\hat P=U^\top Q^{-1/2}PQ^{-1/2}U$, where $U$ is an orthogonal matrix such that         
	$
  	\hat P= \mathrm{diag}( P_1,0 ),
 	$
	with $P_1\succ 0$. Hence, the \blue{matrix inequality} in the statement can be rewritten as  
	$$
	\hat A^\top \hat P +\hat P\hat A+\hat A^\top\hat A \preccurlyeq 0 \,,
	$$        
	Letting	$
	\hat  A= \left [
	\begin{smallmatrix} 
	\hat A_1 & &\hat A_3 \\ 
	\hat A_2 & & \hat A_4 
	\end{smallmatrix} \right ],
	$
the previous inequality implies that $\hat A_3=0$, $\hat A_4=0$, so that 
	$
	\hat  A= \left [
	\begin{smallmatrix} 
	\hat A_1 & & 0 \\ 
	\hat A_2 & & 0
	\end{smallmatrix} \right ]
	$. 
Moreover,  
	$$
	\blue{P_1 \succcurlyeq  \hat A_1^\top P_1 \hat A_1+  W},
	$$
where  
$
\hat  W \coloneqq \hat A_1^\top\hat A_1 + \hat A_2^\top \hat A_2.
$
Due to the assumption on the rank of $P$, $\hat W$ is positive definite, hence $A_1$ is Hurwitz and the system is convergent.
}\hfill$\blacksquare$


\balance

\bibliographystyle{plain}
\bibliography{conv_aut}

\end{document}